\documentclass[10pt]{article}
\usepackage[english]{babel}
\usepackage[latin1]{inputenc}
\usepackage{amscd}
\usepackage{amsfonts}
\usepackage{amsgen}
\usepackage{amsmath}
\usepackage{amssymb}
\usepackage{amstext}
\usepackage{amsthm}
\usepackage{graphicx}
\usepackage{indentfirst}
\usepackage{latexsym}
\usepackage{makeidx}
\usepackage{mathrsfs}
\usepackage{epsfig}
\usepackage{wrapfig}
\usepackage[all,knot,arc,import,poly]{xy}
\usepackage[usenames]{color}
\newtheorem{Df}{Definition}[section]
\newtheorem{Teo}[Df]{Theorem}
\newtheorem{Prop}[Df]{Proposition}

\newtheorem{Ex}[Df]{Example}
\newtheorem{Obs}[Df]{Remark}

\newtheorem{Que}[Df]{Questão}
\newtheorem{Cor}[Df]{Corollary}

\newcommand{\n}{\noindent}
\newcommand{\Dem}{\n{\bf Proof:\;\;}}

\newcommand{\bc}{\begin{center}}
\newcommand{\ec}{\end{center}}

\newcommand{\N}{\mathbb{N}}

\newcommand{\vsete}{\vspace*{0.7cm}}
\newcommand{\vcinco}{\vspace*{0.5cm}}
\newcommand{\vtres}{\vspace*{0.3cm}}

\newlength{\dede}     % bloco de comandos para tabulaca~o

\newcommand{\cF}{\mbox{$\cal F$}}

\newcommand{\cL}{\mbox{$\cal L$}}

\newcommand{\cA}{\mbox{$\cal A$}}

\newcommand{\cQ}{\mbox{$\cal Q$}}

\newcommand{\cS}{\mbox{$\cal S$}}

\def \rest {\restriction}
\def \ra {\rightarrow}

\def \thra {\twoheadrightarrow }

\newcommand{\Sf}{\cS_f}

\newcommand{\Lf}{\cL_f}
\newcommand{\Qf}{\cQ_f}

\newcommand{\Qc}{\cQ_f^c}
\def \beq { \begin{equation} }
\def \eeq { \end{equation} }

\def \rest {\restriction}

\newcommand{\baf}{\begin{Afi}\nf{\sl . }}
\newcommand{\eaf}{\end{Afi}}

\newcommand{\bdf}{\begin{Df}\nf{\bf .}}
\newcommand{\edf}{\end{Df}}
\newcommand{\bte}{\begin{Th}\nf{\bf .}}
\newcommand{\ete}{\end{Th}}  %\eth ja esta definido na table 8.20 (amssymb)
\newcommand{\bco}{\begin{Co}\nf{\bf .}}
\newcommand{\eco}{\end{Co}}
\newcommand{\ble}{\begin{Le}\nf{\bf .}}
\newcommand{\ele}{\end{Le}}
\newcommand{\bpr}{\begin{Pro}\nf{\bf .}}
\newcommand{\epr}{\end{Pro}}
\newcommand{\bex}{\begin{Exa}\nf{\bf .} \rm}
\newcommand{\eex}{\end{Exa}}
\newcommand{\bre}{\begin{Rem}\nf{\bf .} \rm}
\newcommand{\ere}{\end{Rem}}
\newcommand{\bfa}{\begin{Fa}\nf{\bf .} \sl}
\newcommand{\efa}{\end{Fa}}
\newcommand{\bqt}{\begin{Que}\nf{\bf .}}
\newcommand{\eqt}{\end{Que}}

\newcommand{\bct}{\begin{Ct}\nf{\bf .} \rm}
\newcommand{\ect}{\end{Ct}}
\newcommand{\bxa}{\begin{Exa}\nf{\bf .} \rm}
\newcommand{\exa}{\end{Exa}}

\newcommand{\bu}{$\bullet$\ }
\newcommand{\sig}{\mathbb{S}ig}
\newcommand{\cat}{\mathbb{C}at}
\newcommand{\set}{\mathbb{S}et}
\newcommand{\bd}{\bar{\partial}}
\newcommand{\br}{\bar{\rho}}
\newcommand{\vepsilon}{\varepsilon}
\newcommand{\vphi}{\varphi}

\newcommand{\sub}{\subseteq}

\begin{document}

\title{An abstract approach to Glivenko's theorem} %\\{\small (first draft: full proofs are in the typing process)}}

\author{Darllan Concei\c{c}\~{a}o\ Pinto \thanks{Instituto de Matem\'atica e Estat\'istica, University of S\~ao Paulo, Brazil.\ Emails: darllan\_\! math@hotmail.com, hugomar@ime.usp.br }\\
Hugo Luiz Mariano \thanks{Research supported by FAPESP, under the Thematic Project LOGCONS:
Logical consequence, reasoning and computation (number 2010/51038-0).}
}

\date{}

\maketitle

\begin{abstract}
The aim of this work is to provide a special kind of conservative translation between abstract logics, namely an \textit{abstract Glivenko's theorem}. Firstly we define institutions on the categories of logic, algebraizable logics, and Lindenbaum algebraizable logic. In the sequel,  we introduce the notion os Glivenko's context relating two algebraizable logics (respectively, Lindenbaum algebraizable
logics) and  we prove that for each Glivenko's context can be associated an institutions morphism between the corresponding logical institutions. As a consequence of the existence of such institutions morphisms, we have established abstract versions of Glivenko's theorem between those algebraizable logics (Lindenbaum algebraizable logics), generalizing the results presented in \cite{To}.
In particular, considering the institutions of classical logic and of intuitionistic logic, we build a Glivenko's context and
thus an abstract Glivenko's theorem  that is exactly the traditional Glivenko's theorem. Finally we present a category of algebraizable logic with Glivenko's context as morphisms. We can interpret the results of this work as an evidence of the (virtually unexplored) relevance of institution theory in the study of propositional logic.
\end{abstract}

\section{Introdution}

The methods of combination of logics has been the main motivation to consider categories of logics. This allows one not just unify a choice to represent a logical system, as well as to study relations between logics. Among many kinds of possible morphisms between logics, one of most important is the ``conservative translation", i.e., a relation between logics that inter-translate proves. The classical Glivenko's theorem, proved by Valery Glivenko in 1929 that says one can translate the classical logic into intuitionistic logic by means double-negation of classical formulas, is in a certain way a kind of conservative translation. This work actually concerns to establish a \textit{abstract Glivenko's theorem} between algebraizable logics.

The logical and mathematical device that we have used here is the notion of \emph{Institution}. This notion was introduced for the first time by Goguen and Burstall in \cite{GB}. This concept formalizes the informal notion of logical system into a mathematical object. The main (model-theoretical) characteristic is that an institution contains a satisfaction relation between models and sentences that are coherent  under change of notation: That  motivated us to consider an institution of a logic, i.e., an institution for a propositional logic $l$ represents all logic $l'$ such that is {\em equipollent} with $l$ (\cite{CG}).  We introduce, in the subsequent sections,  institutions for  abstract logics,  algebraizable logics and  Lindenbaum algebraizable logics.

Concerning the latter, we present the definition of a {\em Glivenko's context} between two algebraizable logics. Recalling , we prove in \ref{MorGlivContInsAL} (\ref{MorGlivContInsLAL}) that for each Glivenko's context relating  two algebraizable logics (respectively, Lindenbaum algebraizable logics), can be associated a  institutions morphism between the corresponding logical institutions . Moreover, in \ref{AbstGliInsAL} (\ref{AbstGliInsLAL}) we have that a Glivenko's context between institutions of algebraizable logics (Lindenbaum algebraizable logics) provides an abstract Glivenko's theorem between those logics, generalizing the results presented in \cite{To}. In particular, considering the institutions of classical logic and of intuitionistic logic \ref{traditionalGliv}, we build a Glivenko's context and thus an abstract Glivenko's theorem such that is exactly the traditional Glivenko's theorem.

In the end of this paper we give a brief discussion about the category of algebraizable logics with Glivenko's context as its morphisms.

\section{Preliminaries}

The several process of combining logics were the main motivations for the systematic study of categories of logics. Here the objects are signature and consequence operator pairs, the morphisms are translations between logics.

\begin{Df}
A signature is a sequence of pairwise disjoint sets $\Sigma=(\Sigma_{n})_{n\in \N}$. In what follows, $X=\{x_{0},x_{1},...,x_{n},...\}$ will denote a fixed enumerable set (written in a fixed order). Denote $F(\Sigma)$ (or $Fm$) (respectively $F(\Sigma)[n]$ or $Fm[n]$), the set of $\Sigma$-formulas over $X$ (respec. with exactly $\{x_{0},...,x_{n-1}\}$ variables). In this sense we have that $\bigcup_{k\leq n}Fm[k]$ is the set of formulas such that its set of variable is included in $\{x_{0},...,x_{n}\}$.

A Tarskian consequence relation is a relation $\vdash\subseteq\wp(F(\Sigma))\times F(\Sigma)$, on a signature $\Sigma=(\Sigma_{n})_{n\in \N}$, such that, for every set of formulas $\Gamma,\Delta$ and every formula $\varphi,\psi$ of $F(\Sigma)$, it satisfies the following conditions:
\begin{itemize}
\item[$\circ$]$\bf{Reflexivity}: $If $\varphi\in\Gamma,\ \Gamma\vdash\varphi$
\item[$\circ$]$\bf{Cut}: $If $\Gamma\vdash\varphi$ and for every $\psi\in\Gamma,\ \Delta\vdash\psi$, then $\Delta\vdash\varphi$
\item[$\circ$]$\bf{Monotonicity}: $If $\Gamma\subseteq\Delta$ and $\Gamma\vdash\varphi$, then $\Delta\vdash\varphi$
\item[$\circ$]$\bf{Finitarity}: $If $\Gamma\vdash\varphi$, then there is a finite subset $\Delta$ of $\Gamma$ such that $\Delta\vdash\varphi$.
\item[$\circ$]$\bf{Structurality}: $If $\Gamma\vdash\varphi$ and $\sigma$ is a substitution, then $\sigma[\Gamma]\vdash\sigma(\varphi)$
\end{itemize}
\end{Df}

\vtres

The notion of logic that we consider is:

\begin{Df}
A logic of type $\Sigma$, or a $\Sigma-logic$, is a pair $(\Sigma,\vdash)$ where $\Sigma$ is a signature and $\vdash$ is a Tarskian consequence relation.
\end{Df}

\begin{Df}
\begin{enumerate}
\item Let $L$ be a lattice. A element $a\in L$ is compact if for every directed subset $\{d_{i}\}$ of $L$ we have $a\leq \bigvee_{i}d_{i}\ \Leftrightarrow\ \exists i(a\leq d_{i})$. $L$ is said algebraic if it is complete lattice such that every element is join of compact elements. We denote the category of algebraic lattice by $AL$
\item Let $l=(\Sigma,\vdash)$ be a logic and $A\in \Sigma-Str$. A subset $F$ of $A$ is a $l$-filter is for every $\Gamma\cup\{\varphi\}\subseteq F(\Sigma)$ such that $\Gamma\vdash \varphi$ and every valuation $v:F(\Sigma)\to A$, if $v[\Gamma]\subseteq F$ then $v(\varphi)\in F$. The pair $\langle M,F\rangle$ is then said to be a matrix model of $l$. The set of all matrix model of $l$ is denoted by $Matr_{l}$.
\end{enumerate}
\end{Df}

\begin{Obs}

\begin{enumerate}\label{functorFi}
\item Let $l=(\Sigma,\vdash)$ be a logic. We have the map $Fi_{l}:\Sigma-Str\to Lat$ (where $Lat$ is the category of lattices) such that for any algebra $A$, $Fi_{l}(A)$ is the lattice of all $l$-filters of $A$. Moreover, one can restrict the codomain to the category $AL$ where the compact elements are finitely generated filters. Thus $Fi_{l}$ is a contravariant functor from the category $\Sigma-Str$ to the category $AL$ where given $f\in hom_{\Sigma-Str}(A,B)$, $Fi_{l}(f)=f^{-1}$ (inverse image).
\item Let {\bf K} be a quasivariety. We have the functor $Co_{\bf K}:\Sigma-Str\to AL$ such that for every algebra $A$, $Co_{\bf K}(A)$ is the lattice of all relative congruence of $A$, i.e., the lattice such that the elements are congruences $\theta$ such that $A/\theta\in \bf K$.
\end{enumerate}
\end{Obs}

\subsection{Categories of signatures and logics with flexible morphisms}

We provide here a definition of category of logics. The ideas behind it come from \cite{JKE} \cite{FC}, \cite{BCC1}, \cite{BCC2} and \cite{CG}.

First of all we define the category of signature with flexible morphism $\Sf$. Before to define this category, let us introduce the following notation:

If $\Sigma=(\Sigma_{n})_{n\in \N}$ is a signature, then $T(\Sigma):=(F(\Sigma)[n])_{n\in \N}$ is a signature too.

A {\bf flexible} morphism $f : \Sigma\to\Sigma'$ is a sequence of functions  $f_n^\sharp : \Sigma_n \to F(\Sigma')[n], {n\in\omega}$.

For each signature $\Sigma$ and $n\in \N$, consider the particular flexible morphism:

$
\begin{array}{lll}
(j_{\Sigma})_{n}:&\Sigma_{n}\to& F(\Sigma)[n]\\
&c_{n}\mapsto&c_{n}(x_{0},...,x_{n-1})
\end{array}
$

For each flexible morphism $f:\Sigma\to\Sigma'$, there is only one function $\check{f}:F(\Sigma)\to F(\Sigma')$, called the extension of $f$, such that:
\begin{itemize}
\item[(i)]$\check{f}(x)=x$, if $x\in X$;
\item[(ii)]$\check{f}(c_{n}(\psi_{0},...,\psi_{n-1}))=f(c_{n})(x_{0},...,x_{n-1})[x_{0}|\check{f}(\psi_{0}),...,x_{n-1}|\check{f}(\psi_{n-1})]$, if $c_{n}\in \Sigma_{n}, n \in \N$.
\end{itemize}

\begin{Df}
The category $\Sf$ is the category of signatures and $flexible$ morphism as above. The composition in $\Sf$ is given by $(f'\bullet f'')^{\sharp}:=((\check{f}\rest\circ f^{\sharp})_{n})_{n\in\omega}$. The identity $id_{\Sigma}$ in $\Sf$ is given by $(id_{\Sigma})^{\sharp}_{n}:=((j_{\Sigma})_{n})_{n\in\omega}$
\end{Df}

\begin{Df}
If $l=(\Sigma,\vdash),l'=(\Sigma',\vdash')$ are logics then a flexible translation morphism $f:l\to l'$ is a flexible signature morphism $f :\Sigma\to\Sigma'$ in $\Sf$ such that ``preserves the consequence relation", that is, for all $\Gamma\cup\{\psi\}\subseteq F(\Sigma)$, if $\Gamma\vdash\psi$ then $\check{f}[\Gamma]\vdash'\check{f}(\psi)$.

The category $\Lf$ is the category of propositional logics and flexible translations as morphisms. Composition and identities are inherent from $\Sf$.
\end{Df}

\subsection{Other categories of logics}

Due to some difficult that was found in the categories of logics mentioned above, are presented in \cite{MaMe}  others categories of logics that  overcome these ``defects".

\begin{Obs} \label{othercat}

\begin{enumerate}

\item[(I)] Still on the category $\Lf$ we have the ``congruential" \{also called {\em selfextensional}\} logics $\Lf^{c}$. This category is a subcategory of $\Lf$ where the logics satisfy the congruence property, i.e., logics that satisfies: \[\varphi_{0}\dashv\vdash\psi_{0},...,\varphi_{n-1}\dashv\vdash\psi_{n-1}\Rightarrow c_{n}(\varphi_{0},...,\varphi_{n-1})\dashv\vdash c_{n}(\psi_{0},...,\psi_{n-1}).\]

    The inclusion functor $\Lf^{c}\hookrightarrow\Lf$ has a left adjoint given by congruential closure operator $l \mapsto l^{(c)}$, i.e., given a logic $l\in\Lf$, $l^{c}$ is the least logic with the same signature of $l$ but the Tarkian relation is congruential.

A morphism  $f : l \ra l' \in \Lf$ is called {\em dense}, when
$\forall \varphi'_n \in F(\Sigma')[n]$ \ $\exists \varphi_n \in F(\Sigma)[n]$ such that\
 $\varphi'_n \dashv'\vdash \check{f}(\varphi_n)$. If $l' \in \mathcal{L}^{c}_{f}$, then
 $f$ is dense iff $\forall c'_n \in \Sigma'_n$ \ $\exists \varphi_n \in F(\Sigma)[n]$ such that\
  $c'_n(x_0, \ldots, x_{n-1}) \dashv'\vdash \check{f}(\varphi_n)$.

\item[(II)] On the category $\Lf$, consider $Q\Lf$  the quotient category by the congruence relation\footnote{I.e., this category has the same class of objects that $\Lf$, and an arrow between $l \to l'$ the logics is an equivalence class of $\Lf$-arrows $f : l \to l'$.}: $f,g\in \Lf(l,l')$, $f\sim g\ iff\ \check{f}(\varphi)\dashv'\vdash\check{g}(\varphi)$. Thus,  by Proposition 4.3 in \cite{CG}, two logics $l,l'$ are equipollent if only if $l$ and $l'$ are $Q\Lf-$isomorphic. All presentation of classical logic are $Q\Lf-$isomorphic.

\item[(III)] In \cite{MaMe} we found the category $Q\Lf^{c}$ quotient of $\Lf^{c}$ (or simply $\Qf^{c}$).

 For $h \in \mathcal{L}^{c}_{f}(l,l')$,  $[h] \in \Qc(l,l')$ is $\Qc$-isomorphism iff   $h$ is a dense morphism and $h$ is a conservative translation\footnote{I.e., $\Gamma \vdash \psi \ \Leftrightarrow \ \check{h}[\Gamma] \vdash' \check{h}(\psi)$, for all $\Gamma \cup\{\psi\} \subseteq F(\Sigma)$.}.

 This category of logics satisfies {\em simultaneously} certain natural conditions:\\
{\em (a)} \ it is canonically related to the major part of logical systems; \\
{\em (b)} \ has good categorial properties (e.g., it is complete, cocomplete and accessible categories); \\
{\em (c)} \ allow a natural notion of algebraizable logical system (\cite{BP},\cite{Cze1}); \\
{\em (d)} \ allow satisfactory treatment of the ``identity problem" of logics.\\

\end{enumerate}

\end{Obs}

\subsection{Categories of algebrizable logics}

%The idea behind of algebraizing a logic emerged trying to connect two independent approach to logic, the equivalence of logic and the assertion and inference. On the two approaches and the ideas of Hilbert, begun the attempts of connect them.

Traditionally algebraic logic has focused on the algebraic investigation of particular classes of algebras related, in some way, to   logics, whether or not they could be connected to some known assertional system by means of the Lindenbaum-Tarski method. However, when such a connection could be established,
there was interest in investigating the relationship between various meta-logical properties of the logical system and the algebraic properties of the
associated class of algebras.

The  Lindenbaum-Tarski method of algebrization of a logic,  associate with the logic a convenient  quotient of  the formula algebra of the logic, by the congruence relation of interprovability: this idea works in classical logic and in some systems of intuitionistic and modal logics.  However this method cannot algebraize other logics. Thus in the end of the 1980's, Blok-Pigozzi  (\cite{BP}) provide a general definition that, in some sense, encompass the traditional method.

Henceforth ``algebraizable logic" will mean ``algebraizable logic in the Blok-Pigozzi sense".

\begin{Df}
Let $\Sigma$ be a signature. We will denote by $\Sigma-Str$ the category with objects given by all the structures (or algebras) on the signature $\Sigma$ and morphisms $\Sigma$-homomorphisms between them. A fundamental example of $\Sigma$-structure is  $F(\Sigma)$, the absolutely free $\Sigma$-algebra on the set of variables $X$.
\end{Df}

\begin{Df}
Given a class of algebras $\textbf{K}$ over the signature $\Sigma$,  the equational consequence associated with $\textbf{K}$ is the relation $\models_{\textbf{K}}$ between a set of equations $\Gamma$ and a single equation $\varphi\equiv \psi$  over $\Sigma$ defined by:
\[\Gamma\models_{\textbf{K}}\varphi\equiv \psi\ i\!f\!\!f\ for\ every\ A\in\textbf{K}\ and\ every\ \Sigma-homomorphism \ \ h:F(\Sigma)\to A,\]
\[\ if \ h(\eta)=h(\nu)\ for\ all\ \eta\equiv\nu\in \Gamma,\ then\ h(\varphi)=h(\psi).\]
\end{Df}

\begin{Df} \label{BP-def}
Let $l=(\Sigma,\vdash)$ be a logic and $\textbf{K}$ be a class of $\Sigma-$algebra. $\textbf{K}$ is a equivalent algebraic semantics for $l$ if $\vdash$ can be faithfully interpreted in $\models_{\textbf{K}}$ in the following sense:
\begin{enumerate}
\item[(1)]there is a finite set $\tau(p)=\{(\delta_{i}(p),\epsilon_{i}(p)),i=1,...,n\}$ of equations in a single variable $p$ such that for all $\Gamma\cup\{\varphi\}\subseteq F(\Sigma)$ and for $j<n$ has been:

    $\Gamma\vdash\varphi\Leftrightarrow\{\tau(\gamma):\gamma\in\Gamma\}\models_{\textbf{K}}\tau(\varphi)$ where $\tau(\varphi)=\{(\delta_{i}(p)[p/\varphi],\epsilon_{i}(p)[p/\varphi]),i=1,...,n\}$.
\item[(2)]there is a finite system $\Delta_{j}(p,q),j=1,...,m$ of two variables formulas (formed by derived binary connectives) such that for all equation $\varphi\equiv\psi,$\\
$\varphi\equiv\psi=\mid_{\textbf{K}}\models\tau(\varphi\Delta\psi)$\\
where $\varphi\Delta\psi=\Delta(\varphi,\psi)$,  $\Delta(\varphi,\psi)=\{\Delta_{j}(\varphi,\psi),j=1,...,m\}$ and $\tau(\varphi\Delta\psi)=\{\delta_{i}(\Delta_{j}(\varphi,\psi)) \equiv \epsilon_{i}(\Delta_{j}(\varphi,\psi));\ i=1,...,n\ and\ j=1,...,m\}$.
\end{enumerate}
In this case we shall say that a logic $l$ is algebraizable. The set $\langle\tau(p),\Delta(p,q)\rangle$ (or just $\langle\tau,\Delta\rangle$) is called an ``algebraizing pair", with $\tau = (\delta, \epsilon)$ as the ``defining equations" and $\Delta$ as the ``equivalence formulas".
\end{Df}

\begin{Prop}\label{BP-def-alter}
Let $\bf K$ an equivalent algebraic semantic for the algebrizable logic $a=(\Sigma,\vdash)$ with algebraizing pair $\langle\tau,\Delta\rangle$, then:
\begin{enumerate}
\item[1.]For all set of equations $\Gamma$ and for all equation $\varphi\equiv\psi$, we have that \[\Gamma\models_{\bf K}\varphi\equiv\psi\ \Leftrightarrow\ \{\xi\Delta\eta:\xi\equiv\eta\in\Gamma\}\vdash\varphi\Delta\psi\]
\item[2.]For each $\psi\in F(\Sigma)$ we have that
\[\psi\dashv \ \vdash\Delta(\tau(\psi)).\]
\end{enumerate}
Conversely, if there is a logic $a=(\Sigma,\vdash)$ and formulas $\langle\Delta(p,q),\tau(p)\rangle$ that satisfy the conditions 1. and 2., then $\bf K$ is an equivalent algebraic semantics for $a$.
\end{Prop}

\begin{Obs} \label{detach} By a direct application of the definition above, if $l=(\Sigma,\vdash)$ is an algebraizable logic and $\phi, \psi \in F(\Sigma)$, then $\phi, \phi \Delta \psi \vdash \psi$ (detachment property).
\end{Obs}

As examples of algebraizable logics we have, in addition to CPC (Classic Propositional Calculus) and IPC (Intuitionistic Propositional Calculus), some modal logics, the Post and Lukasiewicz multi-valued logics, and many of several versions of quantum logic.

In case of CPC (respectively IPC), a possible algebraizing pair $\langle\Delta(p,q),\tau(p) \rangle =
\langle\Delta(p,q), (\delta(p),\epsilon(p))\rangle$ is:

\begin{enumerate}
\item[1.]$\Delta(p,q)=\{p\leftrightarrow q\}$
\item[2.]$\epsilon(p)=p$
\item[3.]$\delta(p)=\top$
\end{enumerate}

and {\bf K} is the class of Boolean algebras (respectively the class of Heyting algebras). %Due to Fact \ref{BP-def -alter}  we have that the two different algebrazing pair for CPC and IPC given above create a same class of algebra for its respective logics.

Recall that a quasivariety is a class of algebras $K$ such that it is axiomatizable by quasi-identities, i.e., formulas of the form
\[(p_{1}\equiv q_{1}\wedge...\wedge p_{n}\equiv q_{n})\to p\equiv q\ for\ n\geq 1\]
when $n=0$ the quasi-identity is \[\top\to p\equiv q.\footnote{That is equivalent to the equation $p\equiv q$.}\]

Now we will recall a result about ``uniqueness" of algebraizing pair and the quasivariety semantics  of an algebraizable logic. For any class $K$ of $\Sigma$-algebras let us denote $(K)^{Q}$ the $\Sigma$-quasivariety generated by $K$.

\begin{Prop} [2.15-\cite{BP}] \label{uniq} Let $a$ be an algebraizable logic.

(a) Let $\langle (\delta_i(p), \varepsilon_i(p)),\Delta_i(p,q)\rangle$, an algebraizing pair for $a$, and $K_{i}$ an equivalent algebraic semantic associated with $a$, for each $i\in \{0,1\}$. Then $(K_{0})^{Q}, (K_{1})^{Q}$ are equivalent algebraic semantics for $a$. Moreover, some uniqueness conditions hold: \\
$\bullet$ \ on quasivariety semantics: $(K_{0})^{Q} = (K_{1})^{Q}$;\\
$\bullet$ \  on equivalence formulas: $\Delta_{0}(p,q)\dashv\vdash\Delta_{1}(p,q)$;\\
 $\bullet$ \ on defining equations: $(\delta_0(p) \equiv  \varepsilon_0(p)) \ =\mid_{K}\models \ (\delta_1(p) \equiv  \varepsilon_1(p))$ (where $K := (K_{0})^{Q} = (K_{1})^{Q}$).

(b) Let $\langle (\delta_i(p), \varepsilon_i(p)),\Delta_i(p,q)\rangle$. Suppose that the following conditions holds:\\
$\bullet$ \  $(\delta_0(p), \varepsilon_0(p)),\Delta_0(p,q)\rangle$ is an algebraizing pair for $a$;\\
$\bullet$ \ $\Delta_{0}(p,q)\dashv\vdash\Delta_{1}(p,q)$;\\
 $\bullet$ \ $(\delta_0(p) \equiv  \varepsilon_0(p)) \ =\mid_{(K_{0})^{Q} }\models \ (\delta_1(p) \equiv  \varepsilon_1(p))$.\\
Then $\langle (\delta_1(p), \varepsilon_1(p)),\Delta_1(p,q)\rangle$is an algebraizing pair for $a$ and $(K_{1})^{Q} = (K_{0})^{Q}$.

\end{Prop}

%With the characterization below is possible show if a determined class of algebras is an equivalent algebraic semantics for a given logic. Hence determine if a quasivariety is the unique equivalence algebraic semantic for a logic.

If $a = (\Sigma, \vdash)$ is an algebraizable logic then, by the Proposition above, we can (and we will) denote by $QV(a)$ the unique quasivariety on the signature $\Sigma$ that is an equivalent algebraic semantics for $a$.

\begin{Prop} [2.17 \cite{BP}] \label{axiomat}
Let  $a$ be an algebraizable logic $a$ and  $\langle (\delta, \epsilon),\Delta\rangle$ be an algebraizing pair for $a$. Then the quasivariety $QV(a)$ is axiomatized by the set given by the 3 kinds of quasi-equations  below:\\
\bu $\delta(x_0 \Delta x_0) \equiv \epsilon(x_0 \Delta x_0)$;\\
\bu $\delta(x_0 \Delta x_1) \equiv \epsilon(x_0 \Delta x_1) \ \to \ x_0 \equiv x_1$;\\
\bu $(\bigwedge_{i < n} \delta(\psi_i) \equiv \epsilon(\psi_i)) \ \to \  \delta(\phi) \equiv \epsilon(\phi)$, for each $\{\phi, \psi_0, \cdots, \psi_{n-1}\} \subseteq F(\Sigma)$ such that $\{ \psi_0, \cdots, \psi_{n-1}\} \vdash \phi$, for $n \geq 0$.
\end{Prop}

An attempt to determine if a given logic is algebraizable, at times found difficulties about the definition given above. Thus we have the following characterization.

\begin{Prop} [4.7-\cite{BP}] \label{conditions}
Let $a=(\Sigma,\vdash)$ be a logic and $\Delta\subseteq_{fin} F(\Sigma)[2]$, $(\delta\equiv\epsilon)\subseteq_{fin} (F(\Sigma)[1]\times F(\Sigma)[1])$ such that the conditions below are satisfied
\begin{enumerate}
\item[(a)]$\vdash\varphi\Delta\varphi$, for all $\varphi\in F(\Sigma)$;
\item[(b)]$\varphi\Delta\psi\vdash\psi\Delta\varphi$, for all $\varphi,\psi\in F(\Sigma)$;
\item[(c)]$\varphi\Delta\psi,\psi\Delta\vartheta\vdash\varphi\Delta\vartheta$, for all $\varphi,\psi,\vartheta\in F(\Sigma)$;
\item[(d)]$\varphi_{0}\Delta\psi_{0},...,\varphi_{n-1}\Delta\psi_{n-1}\vdash c_{n}(\varphi_{0},..., \varphi_{n-1})\Delta c_{n}(\psi_{0},...,\psi_{n-1})$, for all $c_{n}\in\Sigma_{n}$ and all $\varphi_{0},\psi_{0},...,\varphi_{n-1},\psi_{n-1}\in F(\Sigma)$;
\item[(e)]$\vartheta\dashv\vdash\Delta(\tau(\vartheta))$, for all $\vartheta\in F(\Sigma)$.
\end{enumerate}
Then $a$ is an algebraizable logic with $\Delta$ as equivalence formulas and $\tau$ as defining equations.

Conversely if $a=(\Sigma,\vdash)$ is a algebrizable logics with algebraizing pair $\langle\Delta(p,q),\tau(p)\rangle$, then the conditions $(a)$ to $(e)$ are satisfied for these formulas.
\end{Prop}

\begin{Obs} \label{strongeralg}
It follows from the characterization above that, if $\vdash_0,  \vdash_1$ are consequence operators over the same signature $\Sigma$, if $l_0 = (\Sigma, \vdash_0)$ is an algebraizable logic with algebraizing pair $\langle\Delta(p,q),\tau(p)\rangle$ and $\vdash_0 \leq  \vdash_1$ (for any $\Gamma\cup\{\varphi\}$,if $\Gamma\vdash_{0}\varphi$ then $\Gamma\vdash_{1}\varphi$), then $l_1 = (\Sigma, \vdash_1)$ is an algebraizable logic and $\langle\Delta(p,q),\tau(p)\rangle$ is an algebraizing pair.
\end{Obs}

\begin{Df} Let $\Sigma$ be a signature, ${A}$ be a $\Sigma$-algebra  and $F\subseteq A$.

(a) Let $\theta$ be a congruence in ${A}$. $\theta$ is said to be compatible with $F$ if, for all $a,b\in A$, if $a\in F$ and $\langle a,b\rangle\in\theta$  then   $b\in F$. Given an algebra $A$ and a subset $F$ of its domain there always exists a greatest congruence of $A$ compatible with $F$. To prove this it is enough to show that the supremum of the set of congruences of $A$ compatible with $F$ in the complete lattice of congruences of $A$ is a congruence compatible with $F$ (Applying the Zorn's lemma in the set of all compatible congruence?).

(b) We will denote by $\Omega^{A}(F)$ the largest congruence of $A$ compatible with $F$. We say that the function $\Omega^{A}$ with domain the set of all subsets of $A$ is called the Leibiniz operator on $A$.
\end{Df}

With the definition of categories of logics given above, it is possible to define the category of algebraizable logics: its morphisms are the translations of algebraizable logics that preserves algebraizing pairs (note that, by Fact \ref{uniq}, this does not depend on particular choice of algebraizing pair of the source and target logics). Other categories of algebraizable logics can be found in \cite{JKE}, \cite{FC}.

\begin{enumerate}
%\item[$\bullet$]$\cA_{s}$ is the category of algebraizable logics with morphism in $\Ls$ such that preserves algebraizing pair. In the sequence of works, \cite{AFLM1}, \cite{AFLM2}, \cite{AFLM3} is proven that the category $\cA_s$ is a relatively complete $\omega$-accessible category \cite{AR}.
\item[$\bullet$]$\cA_{f}$ is the category of algebraizable logics with morphisms in $\Lf$ such that preserves algebraizing pair. $\cA_f$ is a (non full) subcategory of $\Lf$, $\cA_{f}\hookrightarrow\Lf$.
\item[$\bullet$] Besides the category $\cA_f$, we consider also the following categories: \\
- \ $\cA_{f}^{c} : = \cA_{f} \cap \mathcal{L}^{c}_{f}$, the (sub)category of algebraizable and congruential logics;\\
- \ $Q\cA_{f}$, the quotient category of $\cA_{f}$ by the congruence determined by interdemonstrability relation ($\dashv \ \vdash$);\\
- \  $Q\cA_{f}^{c}$, the quotient category of $\cA_{f}^c$.
\item[$\bullet$]The ``Lindenbaum algebraizable" logics are logics $l\in \cA$ such that given formulas $\varphi,\psi\in F(\Sigma)$, $\varphi\dashv\ \vdash\psi \ \Leftrightarrow \ \vdash\varphi\Delta\psi$ (note this does not depend on the particular choice of $\Delta$; the implication $\Leftarrow$ always hold, by \ref{detach}). The class of Lindenbaum algebraizable logics determines a full subcategory of the category of algebraizable logics ($j:Lind(\cA_{f})\hookrightarrow\cA_{f}$). The category $Lind(\cA_{f})$ plays a relevant role in the representation theory of logics (\cite{MaPi1}, \cite{MaPi2}). The inclusion functor $Lind(\cA_{f})\hookrightarrow\cA_{f}$ has a left adjoint functor $L:\cA_{f}\to Lind(\cA_{f})$. In \cite{MaPi2} is proven that $Lind(\cA_f)=\cA^{c}_{f}$.
\end{enumerate}

\begin{Obs}\label{obs1}
Let $l=(\Sigma,\vdash),l'=(\Sigma',\vdash')$. Given a morphism $h\in hom_{\Lf}(l,l')$ we have a functor $h^{\star}:\Sigma'-str\to\Sigma-str$ defined by:

given $M\in\Sigma'-str$,$h^{\star}(M)$ has the same universe of $M$ and $c_{n}^{h^{\star}(M)}=h(c_{n})^{M}$ for every $c_{n}\in \Sigma_{n}$ and every $n\in\omega$. Now given $f\in hom_{\Sigma-str}(M,N)$, $h^{\star}(f):=f$ is a morphism between $h^{\star}(M)$ and $h^{\star}(N)$. This functor commutes over Set, i.e., commute with the forgetful functors. More detail about that can be found in \cite{MaPi2}
\end{Obs}

\subsection{Institutions and their morphisms}

The notion of \emph{Institution} was introduced for the first time by Goguen and Burstall in \cite{GB}. This concept formalizes the informal notion of logical system into a mathematical object. The main (model-theoretical) characteristic is that an institution contains a satisfaction relation between models and sentences that are coherent  under change of notation: That  motivated us to consider an institution of a logic, i.e., an institution for a propositional logic $l$ represents all logic $l'$ such that is {\em equipollent} with $l$ (\cite{CG}).

We start giving the definition of institution with its notion of morphisms (and comorphisms), and consequently its category.

%$\blacktriangleright$ {\em Institutions}
\begin{Df}
An Institution $I=(\mathbb{S}ig,Sen,Mod,\models)$ consists of

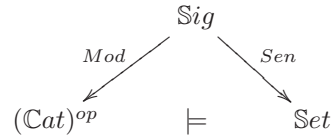
\[\xymatrix{
&\mathbb{S}ig\ar[ld]_{Mod}\ar[rd]^{Sen}&\\
(\mathbb{C}at)^{op}&\models&\mathbb{S}et
}\]

\begin{enumerate}
\item[1.]a category $\mathbb{S}ig$, whose the objects are called {\em signature},
\item[2.]a functor $Sen:\sig\to\set$, for each signature a set whose elements are called {\em sentence} over the signature
\item[3.]a functor $Mod:(\sig)^{op}\to\cat$, for each signature a category whose the objects are called {\em model},
\item[4.]a relation $\models_{\Sigma}\subseteq|Mod(\Sigma)|\times Sen(\Sigma)$ for each $\Sigma\in|\sig|$, called $\Sigma$-{\em satisfaction}, such that for each morphism $h:\Sigma\to\Sigma'$, the compatibility  condition
    \[M'\models_{\Sigma'}Sen(h)(\phi)\ if\ and\ only\ if\ Mod(h)(M')\models_{\Sigma}\phi\]
  holds for each $M'\in|Mod(\Sigma')|$ and $\phi\in Sen(\Sigma)$
\end{enumerate}
\end{Df}

\begin{Ex} \label{inst-examples}
Let $Lang$ denote the  category of  languages $L = ((F_n)_{n \in \N},(R_n)_{n \in \N})$, -- where  $F_n$ is a set of symbols of $n$-ary function symbols and $R_n$ is a   set of symbols of $n$-ary relation symbols, $n \geq 0$ -- and language morphisms\footnote{That can be chosen ``strict" (i.e., $F_n\mapsto F_n'$, $R_n \mapsto R'_n$) or chosen be ``flexible" (i.e., $F_n\mapsto \{n-ary-terms(L')\}$, $R_n \mapsto \{n-ary-atomic-formulas(L')\}$).}. For each pair of cardinals $\aleph_0 \leq \kappa, \lambda \leq \infty$, the category $Lang$  endowed with the usual notion of  $L_{\kappa,\lambda}$-sentences (=  $L_{\kappa,\lambda}$-formulas with no free variable), with the usual association of category of structures and  with the usual (tarskian) notion of satisfaction, gives rise to an institution $I({\kappa,\lambda})$.

\end{Ex}

\begin{Df} Let $I$ and $I'$ be institutions.

(a) An Institution {\bf morphism} $h = (\Phi,\alpha,\beta) : I \to I'$ consists of

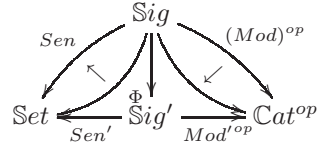
\[\xymatrix{
&\mathbb{S}ig\ar@/^1pc/[ld]_{\nwarrow}\ar@/_/[ld]_{Sen}\ar@/^/[rd]^{(Mod)^{op}}\ar@/_1pc/^{\swarrow}[rd]\ar[d]_>{\Phi}&\\
\set&\sig'\ar[l]^{Sen'}\ar[r]_{{Mod'}^{op}}&\cat^{op}
}\]

\begin{enumerate}
\item[1.]a functor $\Phi:\sig\to\sig'$
\item[2.]a natural transformation $\alpha:Sen'\circ \Phi\Rightarrow Sen$
\item[3.]a natural transformation $\beta:Mod\Rightarrow Mod'\circ\Phi^{op}$
\end{enumerate}
such that the following {\em compatibility condition} holds:
\[m\models_{\Sigma}\alpha_{\Sigma}(\varphi')\ \ if\!f\ \ \beta_{\Sigma}(m)\models'_{\Phi(\Sigma)}\varphi'\]

For any $\Sigma\in\sig$, any $\Sigma$-model $m$ and any $\Phi(\Sigma)$-sentence $\varphi'$.
\vcinco

(b) A triple $f=\langle \phi,\alpha,\beta\rangle:I\to I'$ is a {\bf comorphism} between the given institutions if the following conditions hold:
\begin{itemize}
\item[$\bullet$]$\phi:\sig\to\sig'$ is a functor.
\item[$\bullet$]$\alpha:Sen\Rightarrow Sen'\circ\phi$ and $\beta:Mod'\circ\phi^{op}\Rightarrow Mod$ are natural transformations such that satisfy:

\[m'\models'_{\phi(\Sigma)}\alpha_{\Sigma}(\varphi)\ iff\ \beta_{\Sigma}(m')\models_{\Sigma}\varphi\]

For any $\Sigma\in\sig$,  $m'\in Mod'(\phi(\Sigma))$ and $\varphi\in Sen(\Sigma)$.
\end{itemize}

\end{Df}

\begin{Ex} \label{inst-morph-examples} Given two pairs of cardinals  $(\kappa_i, \lambda_i)$, with  $\aleph_0 \leq \kappa_i, \lambda_i \leq \infty$, $i =0,1$, such that $\kappa_0\leq \kappa_1$ and $\lambda_0  \leq \lambda_1$, then it is induced a morphism and a comorphism of institutions
$(\Phi, \alpha, \beta): I(\kappa_0, \lambda_0) \to I(\kappa_1, \lambda_1)$, given by the {\em same data}: $\sig_{0} = Lang = \sig_{1}$, $Mod_0 = Mod_1 : (Lang)^{op} \to \cat$, $Sen_{i}=L_{\kappa_{i},\lambda_{i}},\ i=0,1$, $\Phi = Id_{Lang} : \sig_0 \to \sig_1$, $\beta := Id : Mod_{i} \Rightarrow Mod_{1-i}$,   $\alpha := inclusion : Sen_{0} \Rightarrow Sen_{1}$.

\end{Ex}

\section{Institutions for abstract propositional logics}

In the first subsection of this section, we provide a \emph{institution} for a  category of propositional logics. That is naturally interesting because the theory of institutions was firstly used by computer scientist for first order logic.

However, the main motivation for  the use of institution theory in this work is because it relates the sentences and models of a logic {\em independently of its presentations}, retaining only its ``essence". More precisely, in the second subsection,  we are going to define  institutions for each (equivalence class of) algebraizable logic and Lindenbaum algebraizable logic: this will enable us to  apply notions and results from institutions to study meta-logic properties of a (equivalence class of) well-behaved logic, as we will exemplify in the next section.
\vtres

\subsection{An institution for the abstract propositional logics}

%Given a logic $l=(\Sigma,\vdash)$ consider the following:
From  to the category of logics ${\cal L}_f$, we define:
\vtres

$\bullet$ $\sig_{} := {\cal L}_f$, the category of propositional logics $l=(\Sigma,\vdash)$ and flexible morphisms.
%is the category which the objects $|\sig_{}|$ is composed by logics as $l'=(\Sigma',\vdash')$ such that is isomorphic with $l$ in the category $Q\Lf$. A morphism $f\in Mor_{\sig_{l}}(l_{1},l_{2})$ is also a isomorphism in $Q\Lf$.
\vtres

$\bullet$ $Sen_{}:\sig_{}\to\set$ where $Sen_{}(l)=\mathcal{P}(F(\Sigma))\times F(\Sigma)$ and given $f\in Mor_{\sig}(l_{1},l_{2})$ then $Sen_{}(f):Sen_{}(l_{1})\to Sen_{}(l_{2})$ is such that $Sen_{}(f)(\langle\Gamma,\varphi\rangle)=\langle\check{f}[\Gamma],\check{f}(\varphi)\rangle$. It is easy to see that $Sen_{}$ is a functor.
\vtres

$\bullet$ $Mod_{}:\sig_{}\to\cat^{op}$ where $Mod_{}(l)=Matr_{l}$ and given $f\in Mor_{\sig}(l_{1},l_{2})$, $Mod_{}(f):Matr_{l_{2}}\to Matr_{l_{1}}$ such that $Mod_{}(f)(\langle M,F\rangle)=\langle f^{\star}(M),F\rangle$. $Mod_{}(f)$ is well defined, indeed:

It is enough to prove that given $\langle M,F\rangle\in Matr_{l_{2}}$, then $F$ is a $l_{1}$-filter in $f^{\star}(M)$. Let $\Gamma\cup\{\varphi\}\subseteq F(\Sigma_{1})$ such that $\Gamma\vdash_{1}\varphi$. Let $v:F(\Sigma_{1})\to f^{\star}(M)$ and suppose that $v[\Gamma]\subseteq F$. We define $\bar{v}:F(\Sigma_{2})\to M$ where $\bar{v}(x)=v(x)$ for all variable $x$ and $\bar{v}(c_{n}(\psi_{0},...,\psi_{n-1}))=c_{n}^{M}(\bar{v}(\psi_{0}),...,\bar{v}(\psi_{n-1}))$ for all formula $\varphi=c_{n}(\psi_{0},...,\psi_{n-1})$ where $c_{n}$ is a n-ary connective. As we saw in the Chapter 2, the function $\check{f}:F(\Sigma_{1})\to f^{\star}(F(\Sigma_{2}))$ is a morphism in $\Sigma_{1}-Str$. Therefore the following diagram commutes

\[\xymatrix{
F(\Sigma_{1})\ar[r]^{v}\ar[d]_{\check{f}}&f^{\star}(M)\\
f^{\star}(F(\Sigma_{2}))\ar[ur]_{\bar{v}}
}
\]

This follows directly from results in \cite{MaPi1}, since $\check{f} = \eta_{f^\star}(X) : F(\Sigma_{1})(X)\to f^{\star}(F(\Sigma_{2})(X))$ is the unity of the adjunction between $\Sigma_1-Str$ and $\Sigma_2-Str$, described in Chapter 2. Anyway, we provide here a more explicit proof: For any variable $x$ we have that $\bar{v}\circ\check{f}(x)=v(x)$. Now suppose that for a formula $c_{n}(\psi_{0},...,\psi_{n-1})$ we have $\bar{v}\circ \check{f}(\psi_{i})=v(\psi_{i})$ with $i\in\{0,...,n-1\}$ then

\[\begin{array}{rcl}
\bar{v}\circ\check{f}(c_{n}(\psi_{0},...,\psi_{n-1}))&=&\bar{v}(f(c_{n})(\check{f}(\psi_{0}),...,\check{f}(\psi_{n-1})))\\
&=&f(c_{n})^{M}(\bar{v}\circ\check{f}(\psi_{0}),...,\bar{v}\circ\check{f}(\psi_{n-1}))\\
&=&c_{n}^{f^{\star}(M)}(\psi_{0}),...,\bar{v}\circ\check{f}(\psi_{n-1}))\\
&=&c_{n}^{f^{\star}(M)}(v(\psi_{0}),...,v(\psi_{n-1}))\\
&=&v(c_{n}(\psi_{0},...,\psi_{n-1}))
\end{array}
\]

Since $v[\Gamma]\subseteq F$ we have $\bar{v}\circ\check{f}[\Gamma]\subseteq F$. $f$ is a morphism between logics, so $\check{f}[\Gamma]\vdash_{2}\check{f}(\varphi)$. Since $\langle M,F\rangle\in Matr_{2}$ therefore $\bar{v}\circ\check{f}(\varphi)\in F$. Hence $F$ is a filter of $l_{1}$.
\vtres

$\bullet$ Given $l\in\sig_{}$  We define a relation $\models\subseteq|Mod_{}(l)|\times Matr_{l}$ as:

Given $\langle M,F\rangle\in Mod_{}(l)$ and $\langle\Gamma,\varphi\rangle\in Sen_{}(l)$,
\[\langle M,F\rangle\models_{l}\langle\Gamma,\varphi\rangle\ if\!f\ for\ all\ v:F(\Sigma_{l})\to M,\ if\ v[\Gamma]\subseteq F,\ then\ v(\varphi)\in F.\]

Now we prove  that $\models$ satisfies the compatibility condition. Let $f:l_{}\to l'_{}$ be a morphism in $\sig_{}$, $\langle M',F'\rangle\in Mod_{}(l'_{})$ and $\langle\Gamma,\varphi\rangle\in Sen_{}(l_{})$.

The universal property of $\check{f}$ defines a bijection:
\begin{center}
$v' \in \Sigma'-Str(F(\Sigma')(X),M') \leftrightsquigarrow v \in \Sigma-Str(F(\Sigma)(X), f^{\star}(M') )$
\end{center}
such that the diagram of functions below commutes

\[\xymatrix{
f^{\star}(F(\Sigma'_{}))\ar[r]^{v'}&f^{\star}(M')\\
F(\Sigma_{})\ar[u]^{\check{f}}\ar[ur]_{v}
}\]

Thus

\[\begin{array}{rcl}
\langle f^\star(M'),F'\rangle\models_{l}\langle\Gamma,\varphi\rangle\ &if\!f& for\ all\ v:F(\Sigma_{})\to f^\star(M'),\ if\ v[\Gamma]\subseteq F',\ then\ v(\varphi)\in F'\\
&if\!f& for\ all\ v':F(\Sigma'_{})\to M',\ if\ v'[\check{f}[\Gamma]]\subseteq F',\ then\ v'(\check{f}(\varphi))\in F'\\
&if\!f& \langle M',F'\rangle\models_{l'}\langle\check{f}[\Gamma],\check{f}(\varphi)
\end{array}
\]

%$``\Rightarrow''$ Let $v:F(\Sigma_{2})\to M$ such that $v(\check{f}[\Gamma])\subseteq F$. consider $v'=v\circ \check{f}$. We have the following diagram commuting:

%\[\xymatrix{
%f^{\star}(F(\Sigma_{2}))\ar[r]^{v}&f^{\star}(M)\\
%F(\Sigma_{1})\ar[u]^{\check{f}}\ar[ur]_{v'}
%}\]

%$v$ is a morphism in $\Sigma_{1}-Str$ and since $\langle f^{\star}(M),F\rangle\in Matr_{1}$ and $v(\check{f}[\Gamma])\subseteq F$, then $v'[\Gamma]\subseteq F$. As $\langle f^{\star}(M),F\rangle\models_{1}\langle\Gamma,\varphi\rangle$, thus $v(\varphi)\in F$. Therefore $\langle M,F\rangle\models_{2}\langle\check{f}[\Gamma],\check{f}(\varphi)\rangle$.
%\vtres

%$``\Leftarrow''$ Now let $v:F(\Sigma_{1})\to f^{\star}(M)$ such that $v[\Gamma]\subseteq F$. Define the morphism $v':F(\Sigma_{2})\to M$ as $v'(x)=v(x)$ for all variable $x$ and recursively $v'(c_{n}(\psi_{0},...,\psi_{n}))=c_{n}^{M}(v'(\psi_{0}),...,v'(\psi_{n-1}))$. We have the following diagram commuting

%\[\xymatrix{
%F(\Sigma_{1})\ar[d]_{\check{f}}\ar[r]^{v}&f^{\star}(M)\\
%f^{\star}(F(\Sigma_{2}))\ar[ur]_{v'}
%}\]

%Since $v[\Gamma]\subseteq F$ we have $v'\circ\check{f}[\Gamma]\subseteq F$. Since $\langle M,F\rangle\models_{2}\langle\check{f}[\Gamma],\check{f}(\varphi)\rangle$,  thus $v'\circ \check{f}(\varphi)\in F$ and then $v(\varphi)\in F$. Therefore $\langle f^{\star}(M),F\rangle\models_{1}\langle\Gamma,\varphi\rangle$.

\begin{Df}
We denote by $I_{f}=\langle\sig_{},Sen_{},Mod_{},\models\rangle$ the above defined institution of abstract propositional logics  associated with ${\cal L}_f$.
\end{Df}
\vtres

\subsection{(Lindenbaum) algebraizable logics as institutions}

 In this section we define  institutions for each (equivalence class of) algebraizable logic and Lindenbaum algebraizable logic: this will enable us to  apply notions and results from institutions to study meta-logic properties of a (equivalence class of) well-behaved logic, as we will exemplify in the next section.

\subsubsection{The institution of an algebraizable logic}\label{InstitutionAL}

Let $a=(\Sigma,\vdash)$ any algebraizable logic and $\Delta$ any of its a set of equivalence formulas. Given $\varphi\in F(\Sigma)$, consider $\varphi/\Delta$ the class of formulas $\psi$ of $a$ such that $\vdash\varphi\Delta\psi$ (this does not depend on the particular choice of $\Delta$). If $\Gamma \subseteq F(\Sigma)$, still denote $\Gamma/\Delta :=\{\varphi/\Delta;\ \varphi\in\Gamma\}$.
Recall that ${\overline{\cA}}_{f}$ denotes the quotient category of $\cA_f$ by the congruence relation  given by $f, f': a_1\to a_2$, \ $f \equiv f'$ iff for each $\varphi_1 \in F(\Sigma_1)$, $\vdash_2 \check{f}(\varphi_1) \Delta_2 \check{f}'(\varphi_1)$, where $\Delta_2$ is any equivalence formulas for $a_2$ (see  Chapter 1, section 3).
%and $\mathcal{P}(F(\Sigma))\times F(\Sigma)/\delta=\{\langle\Gamma/\Delta,\varphi/\Delta\rangle;\ \langle\Gamma,\varphi\rangle\in\mathcal{P}(F(\Sigma))\times F(\Sigma)\}$

Now fix   $a$  an algebraizable logic. Consider:

$\bullet$ $\sig_{a}$ is the category whose objects are the algebraizable logics isomorphic to $a$ in  ${\overline{\cA}}_{f}$ and the morphisms in $\sig_{a}$  are the isomorphisms in ${\overline{\cA}}_{f}$ (i.e., the equivalence class of$\cA_f$-morphisms $f: a_1\to a_2$ is such that there exists a $\cA_f$-morphism $g:a_2\to a_1$ such that $\vdash_{1}\check{g}\circ\check{f}(\varphi_1)\Delta_{2} \varphi_1$ and $\vdash_{2}\check{f}\circ\check{g}(\psi_2)\Delta_{2}\psi_2$, for each $ \varphi_1 \in F(\Sigma_1), \psi_2 \in F(\Sigma_2)$ ).
\vtres

%Let $a=(\Sigma,\vdash)$ an algebraizable logic and $\Delta$ any of its sets of equivalence formulas. Given $\varphi\in F(\Sigma)$, consider $\varphi/\Delta$ the class of formulas $\psi$ of $a$ such that $\vdash\varphi\Delta\psi$. Denote by $\Gamma/\Delta=\{\varphi/\Delta;\ \varphi\in\Gamma\}$ and define:

$\bullet$ $Sen_{a}:\sig_{a}\to \set$ such that $Sen_{a}(a_1)=\mathcal{P}_{fin}(F(\Sigma_1)/\Delta_1)\times F(\Sigma_1)/\Delta_1$ and given $[h]:a_{1}\to a_{2}$, $sen_{a}([h])(\langle\Gamma/\Delta_{1},\varphi/\Delta_{1}\rangle)=\langle\check{h}[\Gamma]/\Delta_{2},\check{h}(\varphi)/\Delta_{2}\rangle$. This is well defined because, if $h \equiv h'$, then for any $\varphi,\varphi'\in F(\Sigma_1)$, if $\varphi/\Delta=\varphi'/\Delta$ then $\vdash_1\varphi\Delta_1\varphi'$ and since $h,h'$ are represent the same morphism in $\overline{\cA}_f$ we have that $\vdash_2\check{h}(\varphi)\Delta_2\check{h'}(\varphi')$.
\vtres

$\bullet$ $Mod_{a}:\sig_{a}\to\cat^{op}$ is such that $Mod_{a}(a'):=Matr^{*}_{a'}$ and given $[f]:a_{1}\to a_{2}$ we define $Mod_{a}([f]):Matr^{*}_{a_{2}}\to Matr^{*}_{a_{1}}$ where $Mod_{a}([f])(\langle M,F\rangle):=\langle f^{\star}M,F\rangle$, this does not depend on the particular representation of $[f]$. We must prove that $Mod$ is well defined, i.e. that $\langle f^{\star}M,F\rangle$  is a reduced matrix. We saw in the previous subsection that $F$ is a $a_1$-filter for $f^{\star}M$ thus,  firstly,  we prove that $\Omega^{f^{\star}M}(F)$ is a congruence in $M$. Let $(a_{i},b_{i})\in\Omega^{f^{\star}M}(F)$ such that $0\leq i\leq n-1$ and $c_{n}$ a n-ary connective in $a_{2}$; denote $c_n(\vec{x}) := c_n(x_0, \cdots, x_{n-1})$. As $[f]$ is a morphism in $\sig_{a}$, then there exists $g:a_{2}\to a_{1} \in \cA_f$ such that $\vdash_{2}\check{f}\circ \check{g}(c_{n}(\vec{x}))\Delta_{2} c_{n}(\vec{x})$. Since $a_{2}$ is algebraizable logic, we have that $\models_{QV(a_{2})}\check{f}\circ g(c_{n})\approx c_{n}(\vec{x})$. As $\langle M,F\rangle\in Matr^{*}_{a_{2}}$, then $M\in QV(a_{2})$. Hence $g(c_{n})^{f^{\star}M}=\check{f}(g(c_{n}))^{M}=c_{n}^{M}$. We  know that $\Omega^{f^{\star}M}(F)$ is a congruence in $f^{\star}M$, thus $(c_{n}^{M}(a_{0},...,a_{n-1}),c_{n}^{M}(b_{0},...,b_{n-1}))=(g(c_{n})^{f^{\star}M}(a_{0},...,a_{n-1}),g(c_{n})^{f^{\star}M}(b_{0},...,b_{n-1}))\in\Omega^{f^{\star}M}(F)$.
Therefore $\Omega^{f^{\star}M}(F)$ is a congruence on $M$. Moreover, it is compatible with $F$. Hence $\Omega^{f^{\star}M}(F)\subseteq\Omega^{M}(F)=Id_{|M|\times |M|}$. Then $\Omega^{f^{\star}M}(F)=Id_{|f^{\star}M|\times |f^{\star}M|}$, so $\langle f^{\star}M,F\rangle$ is a reduced matrix.
\vtres

$\bullet$ To $\models$ we use here a similar definition as in the  subsection above, namely given $\langle M,F\rangle\in Matr*_{a_1}$ and $\langle\Gamma/\Delta,\varphi/\Delta\rangle\in Sen_{a}(a_1)$ then $\langle M,F\rangle\models \langle\Gamma/\Delta,\varphi/\Delta\rangle$ iff for any valuation $v: F(\Sigma_1)(X)\to M$, if $v[\Gamma]\subseteq F$ then $v(\varphi)\in F$. As $M \in Qv(a_1)$,  this is well-defined, i.e., if $\vdash \theta \Delta \theta'$  then $v(\theta) = v(\theta')$, since  $v$ factors uniquely through the quotient morphism $F(\Sigma_1)(X) \thra F(\Sigma_1)(X)/\Delta$. The proof of the compatibility follows from the same way as in the  subsection above.

\begin{Df}
We denote by $InsAL_{a}=\langle\sig_{a},Sen_{a},Mod _{a},\models\rangle$ the above defined institution. This will be called the algebraizable institution of $a$.
\end{Df}

\vtres

\subsubsection{The institution of a Lindenbaum algebraizable logics}\label{instiutionLAL}

Before define the Institution of Lindenbaum algebraizable logics, we define a notion of satisfiability of class of formulas:

\begin{Df}
Let $a$ be algebraizable logic. Given $M\in QV(a)$, $M\models_{QV(a)}[\varphi]\approx[\psi]$ iff for every valuation $v: F(\Sigma_a)(X)\to M$,
\[v(\varphi')=v(\psi')\ \ such\ that\ \ \varphi'\dashv\vdash\varphi\ \ and\ \ \psi'\dashv\vdash\psi\]
\end{Df}

\begin{Obs}
If $a\in Lind(\cA_{f})$ then, since $F(\Sigma_a)(X)/\dashv\vdash\ = F(\Sigma_a)(X)/\Delta$ is the free $QV(a)$-structure on $X$ (see Remark \ref{freerest}), then $\models_{QV(a)}[\varphi]\approx[\psi]\ \Leftrightarrow\ \models_{QV(a)}\varphi\approx\psi$.
\end{Obs}

Given $a\in Lind(\cA_{f})$. Consider the following maps:

$\bullet$ $\sig'_a$ is the category whose the objects are $a_1=(\Sigma_1,\vdash_1)\in Lind(\cA_{f})$, that are isomorphic to $a$ in the quotient category $QLind(\cA_f) = Q(\cA^c_f)$ and the morphisms are only the isomorphisms in $QLind(\cA_f)$.

$\bullet$ $Mod'_{a}:\sig'^{op}_{a}\to\cat$ such that $Mod'_{a}(a_1)=QV(a_1)$ for all $a_1\in|\sig'_a|$ and $Mod'_{a}(a_1 \overset{[h]}\rightarrow \ a_2) = (QV(a_2) \overset{h^\star\!\!\rest}\rightarrow \ QV(a_1))$ (see \cite{MaPi1}).

$\bullet$ We define now the functor $Sen'_{a}:\sig'_a\to\set$.

Let $a_1 \in|\sig'_a|$. The idea here is to describe a convenient set of tuples that  represents quasi-equations in $\Sigma_1$ (i.e., $Eq_{0}\wedge...\wedge Eq_{n-1}\rightarrow Eq$).

For each $s = ([\varphi_0], \cdots, [\varphi_{n-1}], [\psi])$, a non-empty finite sequence in $F(\Sigma_1)/\dashv \vdash$ (the free $QV(a_1)$-structure on the set $X$) and each $(\tau, \Delta)$, an algebraizable pair of $a_1$,  where $\tau=\{(\varepsilon^{j},\delta^{j});\ j=1,...,m\ for\ some\ m\in\omega\}$, let
$$q(s, (\Delta, \tau)) :=  (([\varepsilon(\varphi_0)], [\delta(\varphi_0)]), \cdots, ([\varepsilon(\varphi_{n-1})],
[\delta(\varphi_{n-1})]), ([\varepsilon(\psi)], [\delta(\psi)]))$$
where the {\em notation} $([\varepsilon(\theta)], [\delta(\theta)])$ abbreviates  the pair of finite sequence of equivalence class of formulas: $([\varepsilon^j(\theta),  [\delta^j(\theta)])_j \ with \ j =1, \cdots, m.$.
Note that, as $a_1$ is a {\em congruential algebraizable  logic}, then: \\
(*) \  If $[\theta] = [\theta']$ (i.e., $\theta \dashv \vdash \theta'$), then $\delta(\theta) \dashv \vdash \delta(\theta')$ and  $\varepsilon(\theta) \dashv \vdash \varepsilon(\theta')$. Thus we have an well defined mapping $\varphi/\Delta \overset{t}\mapsto (\varepsilon(\varphi)/\Delta, \delta(\varphi)/\Delta)$ and $q(s, (\Delta, \tau))$ is well-defined;\\
(**) \ conversely, as $\varphi\dashv\vdash\Delta(\epsilon(\varphi),\delta(\varphi))$, then we have and well defined map $(\varepsilon(\varphi)/\Delta, \delta(\varphi)/\Delta) \overset{r}\mapsto \varphi/\Delta$ and $ r \circ t = id$.

Define $q_s :=  \{ q(s, (\tau, \Delta)) : (\tau, \Delta)$ is an algebraizable pair of $a_1\}$ and  then take $Sen'_{a}(a_1) := \{ q_s : s$  is a non-empty finite sequence in $F(\Sigma_1)/\Delta_1\}$. Note that, by the above remark,  the mapping $s \overset{t} \mapsto q_s$ determine a {\em bijection} between the set of non-empty finite sequences in $F(\Sigma_1)/\Delta$ and $Sen'a(a_1)$

 %is the set all tuples \[q=\{( ([\vec{\alpha}_0], [\vec{\beta}_0]), \cdots,  ([\vec{\alpha}_{n-1}], [\vec{\beta}_{n-1}]);  ([\vec{\alpha}], [\vec{\beta}])), j=1,...m\}\] that represents quasi-equations, i.e., $Eq_{0}\wedge...\wedge Eq_{n-1}\rightarrow Eq$ such that $[\alpha_{i}],[\beta_j]$ belongs to  $F(\Sigma_1)(X)/\dashv\ \vdash$, the free $QV(a_1)$-structure on the set $X$, and $\alpha^{j}_{i} = \varepsilon^{j} (\varphi_i), \beta^{j}_{i} = \delta^{j} (\varphi_i)$, for some algebraizable pair of $a_1$, $(\tau, \Delta)$ whose $\tau=\{(\varepsilon^{j},\delta^{j});\ j=1,...,m\ for\ some\ m\in\omega\}$. Just to simplify notation from now on we will denote the set $q$ by $( ([\alpha_0], [\beta_0]), \cdots,  ([\alpha_{n-1}], [\beta_{n-1}]);  ([\alpha], [\beta]))$ and consequently by $(([\varepsilon(\varphi_{0})],[\delta(\varphi_{0})]),...,([\varepsilon(\varphi_{n-1})]),([\varepsilon(\varphi)],[\delta(\varphi)]))$.

Let $[f]:a_{1}\to a_{2}$ be  an isomorphism in $QLind(\cA_{f})$, in particular $\check{f}/\dashv \vdash :   F(\Sigma_1)/\dashv_1\vdash \to  F(\Sigma_2)/\dashv_2\vdash$ is a bijection. Let $s = ([\varphi_0], \cdots, [\varphi_{n-1}], [\psi])$ be a non-empty finite sequence in $F(\Sigma_1)/\dashv_1\vdash$ and $((\varepsilon,\delta), \Delta)$ be an algebraizable pair of $a_1$. Then  $f*s := ([\check{f}(\varphi_0)], \cdots, [\check{f}(\varphi_{n-1})], [\check{f}(\psi)])$ is a non-empty finite sequence in $F(\Sigma_2)/\dashv_2\vdash$ and the mapping
$$q(s, ((\varepsilon,\delta), \Delta)) \  \mapsto \ q(f*s, ((\check{f}(\varepsilon),\check{f}(\delta)), \check{f}(\Delta)))$$
determines a bijection:  $f^+: q_s \ \overset{\cong}\to\ q_{f*s}$.

 Then $Sen'_{a}([f]): Sen_{a}(a_{1})\to Sen'_{a}(a_{2})$ is given by $Sen'_{a}([f])(q_s) := q_{f*s}$ (this map is well defined). It is straitforward check that $Sen'_{a}:\sig'_a\to\set$ is a functor.

Just to simplify notation, from now on we will denote the any element  of the set $q_s$ by
$( ([\alpha_0], [\beta_0]), \cdots,  ([\alpha_{n-1}], [\beta_{n-1}]),  ([\alpha], [\beta])) = (([\varepsilon(\varphi_{0})],[\delta(\varphi_{0})]),...,([\varepsilon(\varphi_{n-1})]),([\varepsilon(\psi)],[\delta(\psi)]))$.

$\bullet$ Given $a'\in\sig_{a}$, $M'\in QV(a')$ and $q'\in Sen_{a}(a')$, we say that $M'\models^{a}q'$ when, for any (and thus for all!) element $( ([\alpha'_0], [\beta'_0]), \cdots,  ([\alpha'_{n-1}], [\beta'_{n-1}]),  ([\alpha'], [\beta']))$ of $q'$, if
\[M'\models_{QV(a')}[\alpha'_{i}]\approx[\beta'_{i}]\ \forall\ i=0,...,n-1\]
then
\[M'\models_{QV(a)}[\alpha']\approx[\beta']\]

Let $[f]: a_1 \to a_2 \in \sig'_a$, $M_2\in QV(a_2)$ and $q\in Sen'_{a}(a_1)$.  Then, as $[f] :a_1\to a_2$ which is a isomorphism in $QLind(\cA_{f})$, then it is easy to see that
\[M_2\models^{a}Sen(f)(q)\ \Leftrightarrow\ Mod(f)(M_2)\models^{a}q\]

\begin{Df}

Then we have that $InsLAL_{a}=\langle\sig'_{a},Sen'_{a},Mod '_{a},\models'\rangle$ is a institution called the Lindenbaum institution of $a$.
\end{Df}

\begin{Obs} \label{compar- inst}
As can be easily checked, each Lindenbaum algebraizable logic $a$, determines the following {\em comorphism of institutions}: \   $h^a = (\Phi^a,\alpha^a,\beta^a) : InsLAL_{a} \to InsAL_{a}$, where:

$\bullet$ \ $\Phi^a : \sig'_a \to \sig_a$ consists of inclusion of categories: $\Phi^a( a_1 \overset{[h_1]}\to a_2) = a_1 \overset{[h_1]}\to a_2$ ;

$\bullet$ \ $\beta^a : Mod_{a}\circ (\Phi^a)^{op} \Rightarrow Mod'_{a}$, given by, for each $a_1 \in |\sig'_a|$,  $\beta^a(a_1) :  Matr^*_a \to QV(a_1)$ is the forgetful functor;

$\bullet$ $\alpha^a : Sen'_a \Rightarrow Sen_a\circ \Phi^a$, given by, for each $a_1 \in |\sig'_a|$, for each $q \in Sen'_a(a_1)$,  let  $s = ([\varphi_0], \cdots, [\varphi_{n-1}], [\psi])$ be the {\em unique} non-empty finite sequence in $F(\Sigma_1)/\dashv \vdash$ such that $q = q_s$, then $\alpha^a(a_1)(q) := (\{[\varphi_0], \cdots, [\varphi_{n-1}]\}, [\psi]) \in \mathcal{P}_{fin}(F(\Sigma_1)/\Delta_)\times F(\Sigma_1)/\Delta = Sen_{a}(a_1)$.

$\bullet$ It holds the compatiblitity condition: for each $a_1 \in |\sig'_a|$, each $(M,F) \in |Mod_a(\Phi^a(a_1)|= Matr^*(a_1)$  and each $q_s \in Sen'_a(a_1)$

$$(M,F) \models^{I_a} \alpha^a(q_s)\ iff\ M \models^{I'_a} q_s$$

And this follows from:

(+) \ For each $v : X  \to M$ and $\varphi \in F(\Sigma_1)$: \\
\[v(\varphi)\in F \ iff  \ v(\varepsilon(\varphi)) = v(\delta(\varphi)) \footnote{Indeed, as $\varphi\dashv\vdash\Delta(\epsilon(\varphi),\delta(\varphi)),\ then\ v(\varphi)\in F\ iff\ v(\Delta(\epsilon(\varphi),\delta(\varphi)))\in F\ iff\ (v(\epsilon(\varphi)),v(\delta(\varphi)))\in\Omega^{M}(F)=id$.}\]

%PROBLEMA COM COMPATIBILIDADE : sera que para (M,F) matriz reduzida e\phi frmula e v valoraco em M, vale
%v(phi) \in F iff v(\epsilon(phi) =v(delta(phi)

\end{Obs}

\begin{Obs}
One can ask ``why do use different notion of institution of a Lindenbaum algebraizable logic instead of the restrict the notion of institution of algebraizable logic to the class of Lindenbaum algebraizable logic?" The answer to this question is that those institutions seem not be isomorphic, but there are notions of abstract Glivenko's theorem for both of them. This means that we have two different approaches to abstract Glivenko's theorem as follow in the next section. We believe that those two different approaches for the abstract Glivenko's theorem can be applied for special classes of logics, for instance we can use the idea behind of the institution for an algebraizable logic as \ref{InstitutionAL} to provide an institution for an equivalential logic. On the other hand, we can use the idea behind of the institution for a Lindenbaum algebraizable logic as \ref{instiutionLAL} to provide an institution for a truth-equational logic.
\end{Obs}

\section{The abstract Glivenko's theorem}

The  Glivenko's theorem allows one translate the classical logic into the intuitionistic logic by means double negation. More precisely, if  $\Sigma$ be a commom signature for expressing presentations of classical propositional logic (CPC) and intuitionistic propositional logic (IPC) -- for instance, $\Sigma = \{\neg, \to, \wedge, \vee\}$-- and $\Gamma \cup\{\varphi\} \subseteq F(\Sigma)$, then \ $\Gamma\vdash_{CPC} \varphi$ \ iff \ $\neg\neg \Gamma\vdash_{IPC} \neg\neg\varphi$.
Here we generalize the Glivenko's theorem between arbitrary algebraizable logics (Lindenbaum algebraizable logics) using the ideas and notions of the Institution Theory applied to the former defined  institutions for algebraizable logics (Lindenbaum algebraizable logics).
%, without lost the classical result, i.e., when applied to classical and intuitionistic logic we have the classical Glivenko's theorem.

%\begin{Obs}\label{godel}
%Let $a=IPC$ and $a'=CPC$ both Lindenbaum algebraizable logics with the same signature. We have the inclusion morphism $h : IPC\to CPC$. So $h^{\star}_{\rest}:BA\to HA$ has left a adjoint functor $G:HA\to BA$. Observe that $h^{\star}_{\rest}$ is the inclusion functor. Hence given $H\in HA$, $G(H)=H/F_{H}$, where $F_{H}$ is the filter in $H$ generated by the subset $\{a\leftrightarrow \neg\neg a:a\in H\}$.
%Its possible to proof that $G(H) \cong H_{\neg\neg}$, where $H_{\neg\neg}$ denote the poset of regular elements of $H$, that is, those elements $x\in H$ such that $\neg\neg x=x$.

%This fact motivate us to investigate the relation of  G\"odel translation with the left adjunct functor $G$.

%As an application of some of the general results in the present work, we derive in \cite{MaPi2} a generalized "Glinvenko's Theorem" related to an $\cA^c_f$-morphism $h : a \to a'$, whenever holds a simple condition of the unity of the adjunction $(G,h^{\star}_{\rest}) : QV(a') \leftrightarrows QV(a)$.

%\end{Obs}

%ESTABELECER RELACAO COM TORRENS

\vsete

\begin{Obs} \label{left adj}

(a)
Let $a=IPC$ and $a'=CPC$ both Lindenbaum algebraizable logics with the same signature. We have the ``inclusion" morphism $h : IPC\to CPC$. Denote $BA$ and $HA$, the  quasivarieties of Boolean algebras and of Heyting algebras on that commom signature. So $h^{\star}_{\rest} = incl:BA\to HA$ has left a adjoint functor $G:HA\to BA$. Observe that $h^{\star}_{\rest}$ is the inclusion functor. Hence given $H\in HA$, $G(H)=H/F_{H}$, where $F_{H}$ is the filter in $H$ generated by the subset $\{a\leftrightarrow \neg\neg a:a\in H\}$, and the quotient $HA$-homomorphism $q_H : H \thra incl(G(H))$ is the $H$-component of the unity of this adjunction.
It is possible to proof that $G(H) \cong H_{\neg\neg}$, where $H_{\neg\neg}$ denote the (boolean algebra) of regular elements of $H$, that is, those elements $x\in H$ such that $\neg\neg x=x$. Moreover, the surjective $HA$-homomorphism $x \in H \mapsto \neg\neg x \in H_{\neg\neg}$ has $HA$-section $ H_{\neg\neg} \mapsto \neg\neg y \in H$.

%qH:  H.--> H/< x<--> neg neg x>,  com x
%mapsto [x], tem secao [x] \mpato \neg\neg x

(b) Let $h:a\to a'\in \cA_f$. Then $h^{\star}$ and $h^{\star}\!\!\rest$ have respective left adjoints $L_h$ and $\bar{L_h}$. Consider $\partial:Id\Rightarrow h^{\star}\circ L_h$ and $\bar{\partial}:Id\Rightarrow h^{\star}\!\!\rest\circ\bar{L_h}$ the units of the adjunctions between $h^{\star},L_h$ and $h^{\star}\!\!\rest,\bar{L_h}$ respectively. Given $X\in \set$ the following diagram commute: (Here $\partial_{X}=\partial_{FX} = \check{h}$. The same for $\bar{\partial}$)

\[\xymatrix{
FX\ar[r]^{\partial_{X}}\ar[d]&F'X\ar[d]\\
FX/\Delta\ar[r]_{\bar{\partial}_X}&F'X/\Delta'
}\]

Due to results in \cite{MaPi1}, $\partial_{X} = \check{h}$. Moreover, observe that $\bar{\partial}_{X}$ and $[\check{h}]:FX/\Delta\to h^{\star}\!\!\rest(F'X/\Delta')$ both satisfies the universal property, so there exist an isomorphism between $\bar{L}_h(FX/\Delta)$ and $F'X/\Delta$. With this we can consider $\bar{\partial}_{X}:FX/\Delta\to h^{\star}\!\!\rest(F'X/\Delta')$
\end{Obs}

Now we are ready to propose the following

\begin{Df}\label{GlivenkoContext}
A \textbf{Glivenko's context} is a pair $\mathbb{G}=(h:a\to a',\bar{\rho})$ where $h\in \cA_{f}(a,a')$ and $\br : h^{\star}\!\!\rest \circ L_h \Rightarrow Id $ is  a natural transformation that is a {\em section} of the   unit $\bd:Id\Rightarrow h^{\star}\!\!\rest\circ L_h)$.
\end{Df}

%\begin{Ex} Let  $\Sigma$ be a commom signature for expressing presentations of classical propositional logic (CPC) and intuitionistic propositional logic (IPC) -- for instance, $\Sigma = \{\neg, \to, \wedge, \vee\}$) -- and consider  ...by result in cahpter 2
 %then () is a Glivenkio context take $a$ IPC in this languge and $CPC$ and let $h$ the identity map (is a $Lind(\cA_f)$-morphism) that is oviolusy dense regular elemnts

%\end{Ex}

\begin{Obs}\label{GliDense} Let $\mathbb{G}=(h:a\to a',\bar{\rho})$ is a \textbf{Glivenko's context} then:

(a) $[\check{h} = \bar{\partial}_{X} :FX/\Delta\to h^{\star}\!\!\rest(F'X/\Delta')$ is a surjective homomorphism thus $h$ is a $\Delta$-dense morphism (see also cite{MaPi1}). For each $Y \subseteq X$, can be chosen (non naturally) a ``lifting" $\rho_Y : F'Y \to FY$,  for each of the natural sections $\bar{\rho}_Y : F'Y/\Delta_Y \to F'Y/\Delta_Y$:

%\[\xymatrix{
%FX\ar[r]^{\partial_{X}}\ar[d]&F'X\ar[d]\\
%FX/\Delta\ar[r]_{\bar{\partial}_X}&F'X/\Delta'
%}\]

\[\xymatrix{
FY\ar[d]&F'Y\ar[l]_{\rho_{Y}}\ar[d]\\
FY/\Delta &F'Y/\Delta'\ar[l]^{\bar{\rho}_{Y}}
}\]

$\bd_{X}[\theta]=[\check{h}(\theta)]$, for all $\theta \in FX$.

(b)  On the other hand, the condition of being a $\Delta$-dense  on a $\cA_f$-morphism $h$ is not sufficient to ensure that $h$ is part of a Glivenko's context: Consider $a$ the  ``logic of abelian groups" and $a'$ the ``logic of groups" (see Chapter 1, section 3): both are algebraizable logics; then $QV(a) = Ab, QV(a') = Gr$ and, for each group $G$, the unity of this adjunction at $G$ is the quotient homomorphism $q_G : G \thra incl(G/[G,G])$; taking $G = F(x,y)$, the free group in 2 generators, then $G/[G,G] \cong \mathbb{Z} \oplus \mathbb{Z}$ is the free abelian group in 2 generators and is straitforward $q_G : G \thra incl(G/[G,G])$ does not have a section! It will be interesting determine additional   condition on a  $\Delta$-dense morphism, that ensures it  be a part of a Glivenko's context.
%(sao  imcompativeis injtividade e comutatividade em F(x,y):  calcule nos 2 subfgrupos em soma direta dados por n(1,0) e k(0,1) )
% QUESTION: or REMARK : However not every dense morphisms induces the institution morphism $M_{G}$.

(c) Observe that for any $M'\in QV(a')$ there is $M\in QV(a)$ such that
$L_{h}(M)\cong M'$: indeed, as $h : a \to a'$ is a  $\Delta$-dense morphism. We have proved in \cite{MaPi1} that $h^\star\!\!\rest : QV(a') \to QV(a)$ is a full and faithfull functor with a left adjoint and a well-known result on adjunctions, entails that the {\em co-unity} of the adjunction $\kappa$ must be an isomorphism, thus $\kappa_{M'} : L_h(h^\star(M')) \overset{\cong}\to M'$, for each $M' \in QV(a')$.

\end{Obs}

\begin{Obs} \label{idemp-formula}
If $\mathbb{G}=(h:a\to a',\bar{\rho})$ is a \textbf{Glivenko's context} then, taking $Y = \{x_0\} \sub X$, then $E_{Y}(x_0) \in F(Y)$ is a $\Sigma'$-formula in at most one variable $x_0$ such that $[x_0] =  [\check{h}(\rho_{Y}(x_0))] \in F'(Y)/\Delta'$ and thus $[\rho_Y(x_0)] =  [\rho_Y(\check{h}(\rho_{Y}(x_0)))] \in F(Y)/\Delta$.

(Note that the formula $\neg \neg (x)$ appears as a ``fixed formulas" in CPC  and as an ``idempotent formula" in IPC.) Conversely,  give  a ``fixed formula" seems to be also a sufficient condition for exists a Glivenko's context, i.e. give a $\Sigma_{a}$-formula in at most variable $x_0$,  $\theta(x_0)$, such that $ \vdash_a' x_0 \Delta' (\check{h}(\theta(x_0))$. Further investigation is needed to establish (and explore) a precise relation between fixed/idempotent formulas and Glivenko's contexts.

%Neciesaaria ( x  del x0)
%talvez suficiente
%Livre em 1 gerador  livre en  n geradores (comutatriva) (del encolhe var, colimi) quocienbte

\end{Obs}

\begin{Cor} \label{extending Gliv} Let $\mathbb{G}=(h:a\to a',\bar{\rho})$ be a {Glivenko's context} and suppose that $a_1$ is an algebraizable logic and $[e_1] : a \ra a_1$ is an isomorphism in the quotient category $\overline{\cA_f}$. Let $[h_1] : a_1 \ra a'$ be the unique  $\overline{\cA_f}$ such that the diagram below commutes

\[\xymatrix{
a\ar[r]^{[h]}\ar[d]_{[e_1]}&a'\ar[d]^{[id_a']}\\
a_1\ar[r]_{[h_1]}&a'
}\]

Then $h_1 : a_1 \ra a$ is a $\Delta$-dense morphism in $\cA_f$.

From the choice of left adjoints of functors between quasivarieties induced by $\Delta$-dense morphisms  (see Chapter 2), we have the {\em strict} equalities  $L_{h_1} \circ L_{e_1} = L_{h_1 \circ e_1} = L_h$ and then also the diagram below commutes ($L_{e_1}$ is the inverse isomorphism of $e_1^\star$)
\[\xymatrix{
QV(a)\ar[r]^{L_h}&QV(a')\\
QV(a_1)\ar[r]_{L_{h_1}}\ar[u]^{e^\star_1}&QV(a')\ar[u]_ {id_{QV(a')}}
}\]

Thus, the (natural) section, $\bar{\rho}$, of the unity of the adjunction $L_h \dashv h$ induces uniquely a  (natural) section, $\bar{\rho}^{a_1}$, of the unity of the adjunction $L_{h_1} \dashv h_1$.

In more details: if $M_1 \in QV(a_1)$ and $\partial^{a_1}_{M_1} : M_1 \thra h_1^\star(L_{h_1}(M_1))$ is the (canonical) unity of $L_{h_1} \dashv h_1$ (remember that $h_1$ is $\Delta$-dense, since $h$ is $\Delta$-dense and $[e_1]$ is an isomorphism), then
$$e_1^\star(\partial^{a_1}_{M_1}) : e_1^\star(M_1) \thra e_1^\star(h_1^\star(L_{h_1}(M_1))) \ =$$
$$\partial^a_{e_1^\star(M_1)} :  e_1^\star(M_1) \thra h^\star(L_{h}(e_1^\star(M_1)))$$
Thus take $\bar{\rho}^{a_1}_{M_1} := L_{e_1}(\bar{\rho}_{e_1^\star(M_1)})$

\end{Cor}

\subsection{The abstract Glivenko's theorem in \emph{InsAL}}\label{InsAL}

On the category \emph{InsAL} we  are going to present the abstract Glivenko's theorem through morphisms in this category.

\begin{Teo}\label{MorGlivContInsAL}
Let $a, a'$ be algebraizable logics, then each $\mathbb{G}=(h:a\to a',\rho)$ \textbf{Glivenko's context} induces a institutions morphism $ InsAL_{a} \ra InsAL_{a'}$. More precisely, fixing a choice of isomorphisms $\varepsilon: Obj(\sig_{a}) \to Mor(\sig_a)$, $a_1 \mapsto \varepsilon(a_1) = [e_1] : a \overset\cong\to a_1$, we define a institution morphism $N_{(G,\varepsilon_a)} : InsAL_{a}\ra InsAL_{a'}$\footnote{Such induced morphisms are ``isomorphic", for different choices of isomorphisms $\varepsilon^0, \varepsilon^1$.}
\end{Teo}

\Dem

By simplicity, we will write $(G,\varepsilon)$ for $(G,\varepsilon_{a})$. We will define \[N_{(G,\varepsilon)}=\langle \Phi^{(G,\varepsilon)},\alpha^{(G,\varepsilon)},\beta^{(G,\varepsilon)}\rangle\]
(this will depend only on the choice of isomorphisms in the {\em domain} institution $InsAL_{a}$):

$\bullet \ \Phi^{(G,\varepsilon)} : \sig_a \to \sig_{a'}$

The object part of $\Phi^{(G,\varepsilon)}$ is easy do define:  for $a_1 \in |\sig_a|$, set $\Phi^{(G,\varepsilon)}(a_1) :=a'$.

It follows from  adaptations of  results in \cite{AFLM} and  \cite{MaMe} that $\overline{\cA_f}$ is a finitely accessible category that has all colimits (except initial object) and is relatively complete (i.e, has  limits for all diagrams that admits a cone). In particular $\overline{\cA_f}$ has pushouts, and for each $\overline{\cA_f}$-isomorphism $[f] : a \ra a$, we consider the following pushout

\[\xymatrix{
a\ar[r]^{[h]}\ar[d]_{[f]}&a'\ar[d]^{[f'_1]}\\
a\ar[r]_{[h_1]}&a'_1
}\]

%\[\xymatrix{
%a\ar[r]^{h}\ar[d]_{f_1}&a'\ar[d]^{f'_1}\\
%a_1\ar[r]_{h_1}&a'_1
%}\]

As a pushout of an iso is an iso and a pushout of an epi is an epi (recall that $h$ is a $\Delta$-dense morphisms, i.e., $[h]$ is an epi), we may suppose  that the vertex of the pushout is $a'$, $ [f^h] : a' \ra a'$ is an isomorphism and the diagram below commutes\footnote{In this case, this is a necessary and sufficient condition to be a pushout.}

\[\xymatrix{
a\ar[r]^{[h]}\ar[d]_{[f]}&a'\ar[d]^{[f^h]}\\
a\ar[r]_{[h]}&a'
}\]

Note that, as $[h]$ is an epi, then $[f^h]$ is uniquely determined.

Now let $a_1,a_2 \in \sig_a$ and $[g] : a_1 \ra a_2$ be an arrow in $\sig_a$ (i.e., $[g]$ is a $\overline{\cA_f}$-isomorphism). Then, as $ e_i : a \ra a_i$ is an isomorphism, $i =1,2$, then there is a unique isomorphism $[g_\varepsilon] : a \ra a'$ such that left diagram below commutes.

\[\xymatrix{
a_1\ar[d]_{[g]}&a\ar[l]_{[e_{1}]}\ar[r]^{[h]}\ar[d]_{[g_\varepsilon]}&a'\ar[d]^{[g_{\varepsilon}^{h}]}\\
a_2&a\ar[l]^{[e_2]}\ar[r]_{[h]}&a'\\
}\]

\vtres

Then define $\Phi^{(G,\varepsilon)}([g] : a_1 \ra a_2) := [g_\varepsilon^h] : a' \ra a'$. As $[ g_\varepsilon]$ and $[g_\varepsilon^h]$ are {\em uniquely determined} by $g$, it follows that $\Phi^{(G,\varepsilon)}$ preserves identities and composition of arrows in $\sig_a$, thus being a functor.

%\[\xymatrix{
%a\ar[r]^{h}\ar[d]_{f_1}&a'\ar[d]^{f'_1}\\
%a_1\ar[r]_{h_1}&a'_1
%}\]

%Where $f_1$ is a isomorphism between $a$ and $a_1$.

%Observe that this pushout always exist, just take $a'_1=a'$ and $f'_1=Id_{a'}$ (here we use this case).

%The diagram above is a pushout, then $f'_1$ is a isomorphism in $\cQ\cA_{f}$.

%Using the pullback lemma it is easy to show that $\Phi_{h}$ is a functor.

%\[\xymatrix{
%a\ar[r]^{h}\ar[d]_{f_1}&a'\ar[d]^{f'_1}\\
%a_1\ar[r]^{h_1}\ar[d]_{f_2}&a'_1\ar[d]^{f'_2}\\
%a_2\ar[r]^{h_2}&a'_2\\
%}\]
\vtres

$\bullet$ \ $\alpha^{(G,\varepsilon)}: Sen_{a'}\circ\Phi^{(G,\varepsilon)}\Rightarrow Sen_{a}$ where,  for $a$ we have that $\alpha^{(G,\varepsilon)}(a):Sen_{a'}\circ\Phi^{(G,\varepsilon)}(a)=Sen_{a'}(a')\rightarrow Sen_{a}(a)$ such that $\alpha^{(G,\varepsilon)}(a)(\langle\Gamma'/\Delta',\varphi'/\Delta'\rangle)=\langle\rho_{X}[\Gamma']/\Delta,\rho_{X}[\varphi']/\Delta\rangle$. Now for $a_{1}\in\sig_{a}$, let $[e_1] : a \ra a_1$ the isomorphism corresponding by the choice $\varepsilon$ at $a_1$ then,  by \ref{extending Gliv},  $\alpha^{(G,\varepsilon)}(a_1):Sen_{a'}\circ\Phi^{(G,\varepsilon)}({a_1})\to Sen_{a}(a_1)$ such that for $\langle\Gamma'/\Delta',\varphi'/\Delta'\rangle\in Sen_{a'}(a')$, $\alpha^{(G,\varepsilon)}(a_1)(\langle\Gamma'/\Delta',\varphi'/\Delta'\rangle)=\langle\rho^{a_1}_{X}[\Gamma']/\Delta_1,\rho^{a_1}_{X}(\varphi')/\Delta_1\rangle = \langle\check{e_1}\circ\rho_{X}[\Gamma']/\Delta_1,\check{e_1}\circ\rho_{X}(\varphi')/\Delta_1\rangle$. If $[g] : a_1 \ra a_2$ is an isomorphism in $\sig_a$, then for each $\theta' \in F'X$, $\vdash_2 \check{g} (\rho_X^{a_1}(\theta')) \Delta_2 \rho_X^{a_2}(\check{g}^{h}_\varepsilon(\theta'))$, thus  $\alpha^{(G,\varepsilon)}$ is a natural transformation.
%such that $f_{1}:a\to a_{1}$ is a isomorphism in $Q^{c}_{f}$. Notice that for any different choice $f'_{1}:a\to a_1\in Q^{c}_{f}$ we have that $\check{f}_{1}(\varphi)/\Delta_1=\check{f'}_{1}(\varphi)/\Delta_1$ and then we have the well definition of $\alpha_{h}$.  It is easy to see that $\alpha_{h}$ is a natural transformation.
\vtres

%%%%%%%%%%%%%%%  REVER E COMPLETAR   &&&&&&&&&&&&&&&&&&&&&&&&

$\bullet$ \ $\beta^{(G,\varepsilon)}: Mod_{a}\Rightarrow Mod_{a'}\circ(\Phi^{(G,\varepsilon)})^{op}$ where for $a$ we have $\beta^{(G,\varepsilon)}(a):Mod_{a}(a)=Matr^{*}_{a}\to Mod_{a'}(\Phi^{(G,\varepsilon)}(a))=Matr^{*}_{a'}$ such that $\beta_{h}(a)(\langle M,F_{M}\rangle)=\langle L_{h}(M),F_{L_{h}(M)}\rangle$, where $F_{L_{h}(M)} := \bar{\partial}_{M}[F_M]$ (note that $L_{h}(M) \in QV(a)$)\footnote{That $\langle L_{h}(M),F_{L_{h}(M)}\rangle \in Matr^{*}_{a'}$, follows from an argument analogous to the proof of compatibility condition.}. Now for $a_1\in \sig_{a}$, $\beta^{(G,\varepsilon)}(a_1):Mod_{a}(a_1)=Matr^{*}_{a_1}\to Mod_{a'}(\Phi^{(G,\varepsilon)}(a_{1}))=Matr^{*}_{a'}$ such that
$\beta^{(G,\varepsilon)}(a_1)(\langle M,F_{M}\rangle)=\langle L_{h}(e_{1}^{\star}(M')),F_{L_{h}(e_{1}^{\star}(M'))}\rangle$. Similarly of above we have the well definition of $\beta^{(G,\varepsilon)}$. The naturality is proved using the functorial encoding of equipollence that we have proved in Chapter 2.
\vtres

%%%%%%%%%%%%%%%  REVER E COMPLETAR   &&&&&&&&&&&&&&&&&&&&&&&&

$\bullet$ \ The proof the compatibility condition  will be splited in two parts:

{\bf (I)} The first part consist of the compatibility on the logic $a$:

{\bf Claim.} Given $\langle M,F_{M}\rangle\in Mod_{a}(a)=Matr^{*}_{a}$ and $\langle\Gamma'/\Delta',\varphi'/\Delta'\rangle\in Sen_{a'}$ then
\[\beta_{h}(a)(\langle M,F_{M}\rangle)\models'\langle\Gamma'/\Delta',\varphi'/\Delta'\rangle\ \ i\!f\!\!f\ \ \langle M,F_{M}\rangle\models \alpha_{h}(a)(\langle\Gamma'/\Delta',\varphi'/\Delta'\rangle)\]
In other notation
\[\langle L_{h}(M),F_{L_{h}(M)}\rangle\models'\langle\Gamma'/\Delta',\varphi'/\Delta'\rangle\ \ i\!f\!\!f\ \ \langle M,F_{M}\rangle\models\langle \rho_{X}[\Gamma']/\Delta,\rho_{X}(\varphi')/\Delta\rangle\]

{\bf Proof of the Claim.}

$``\Rightarrow"$: Let $v:X\to M$ be an evaluation such that $v[\rho_{X}[\Gamma']]\subseteq F$. We can consider $\bar{w}=h^{*}\circ L_{h}(\bar{v}):(F'(X)/\Delta')^{h}\to h^{*}(L_{h}(M))$ and then the following diagram commutes:

\[\xymatrix{
&X\ar[dl]\ar@/_3pc/[ldd]_{v}\ar[dr]\ar@/^3pc/[rdd]^{w}&\\
\frac{FX}{\Delta}\ar@<1ex>[rr]^{\bd_{X}}\ar[d]^{\bar{v}}&&(\frac{F'X}{\Delta'})^{h}\ar[d]_{\bar{w}}\ar@<1ex>[ll]^{\br_{X}}\\
M\ar@<1ex>[rr]^{\bd_{M}}&&h^{\star}\!\!\rest L_{h}M\ar@<1ex>[ll]^{\br_{M}}
}\]

Since $v\circ \rho_{X}[\Gamma']\subseteq F_{M}$ we have that $\Gamma'\subseteq \rho_{X}^{-1}\circ v^{-1}[F_{M}]$. Consider $(\Delta',\tau')$ a algebraizable pair for $a'$. Then we have that for all $\psi\in\Gamma'$ and $(\varepsilon'^{j},\delta'^{j})\in\tau'$, $(\varepsilon'^{j}(\psi),\delta'^{j}(\psi))\in\Omega^{F'(X)^{h}}(\rho_{X}^{-1}\circ v^{-1}(F_{M}))=\rho_{X}^{-1}\circ v^{-1}(\Omega^{M}(F_{M}))$. Therefore $(v\circ\rho_{X}(\varepsilon'^{j}(\psi)),v\circ\rho_{X}(\delta'^{j}(\psi)))\in\Omega^{M}(F_{M})$. Since $\langle M,F_{M}\rangle$ is a reduced matrix, we have for all $\psi\in\Gamma'$
\[
\begin{array}{rcl}
v\circ\rho_{X}(\varepsilon'^{j}(\psi))&=&v\circ\rho_{X}(\delta'^{j}(\psi))\\
\bar{v}\circ\bar{\rho_{X}}(\varepsilon'^{j}(\psi)/\Delta')&=&\bar{v}\circ\bar{\rho_{X}}(\delta'^{j}(\psi)/\Delta')\\
\bar{\rho_{M}}\circ\bar{w}(\varepsilon'^{j}(\psi)/\Delta')&=&\bar{\rho_{M}}\circ\bar{w}(\delta'^{j}(\psi)/\Delta')\\
\bd_{M}\circ\bar{\rho_{M}}\circ\bar{w}(\varepsilon'^{j}(\psi)/\Delta')&=&\bd_{M}\circ\bar{\rho_{M}}\circ\bar{w}(\delta'^{j}(\psi)/\Delta')\\
\bar{w}(\varepsilon'^{j}(\psi)/\Delta')&=&\bar{w}(\delta'^{j}(\psi)/\Delta')\\
w(\varepsilon'^{j}(\psi))&=&w(\delta'^{j}(\psi))\\
\end{array}\]

Then $(w(\varepsilon'^{j}(\psi)),w(\delta'^{j}(\psi)))\in\Omega^{L_{h}M}(F_{L_{h}M})$. Thus $w(\psi)\in F_{L_{h}M}$ for all $\psi\in\Gamma'$, by assumption $w(\varphi')\in F_{L_{h}M}$. So $(w(\varepsilon'^{j}(\varphi')),w(\delta'^{j}(\varphi')))\in\Omega^{L_{h}M}(F_{L_{h}M})$. Therefore

\[\begin{array}{rcl}
w(\varepsilon'^{j}(\varphi'))&=&w(\delta'^{j}(\varphi'))\\
\bar{w}(\varepsilon'^{j}(\varphi')/\Delta')&=&\bar{w}(\delta'^{j}(\varphi')/\Delta')\\
\bar{\rho_{M}}\circ\bar{w}(\varepsilon'^{j}(\varphi')/\Delta')&=&\bar{\rho_{M}}\circ\bar{w}(\delta'^{j}(\varphi')/\Delta')\\
\bar{v}\circ\bar{\rho_{X}}(\varepsilon'^{j}(\varphi')/\Delta')&=&\bar{v}\circ\bar{\rho_{X}}(\delta'^{j}(\varphi')/\Delta')\\
v\circ\rho_{X}(\varepsilon'^{j}(\varphi'))&=&v\circ\rho_{X}(\delta'^{j}(\varphi'))\\
\end{array}\]

Then $(v\circ\rho_{X}(\varepsilon'^{j}(\varphi')),v\circ\rho_{X}(\delta'^{j}(\varphi')))\in\Omega^{M}(F_{M})$. Therefore $v\circ\rho_{X}(\varphi')\in F_{M}$.
\vtres

$``\Leftarrow"$: Let $w:X\to L_{h}M$ a valuation such that $w[\Gamma']\subseteq F_{L_{h}M}$. Consider $\bar{w}:F'(X)/\Delta'\to L_{h}(M)$ given by $w$ such that the following diagram commutes:

\[\xymatrix{
X\ar[r]\ar[dr]_{w}&F'(X)\ar[r]^{q}\ar[d]^{w}&F'(X)/\Delta\ar[dl]^{\bar{w}}\\
&M&
}\]

Let $\bar{v}=\bar{\rho}_{M}\circ\bar{w}\circ\bd_{X}$, then $\bd_{M}\circ\bar{v}=\bd_{M}\circ\bar{\rho}\circ\bar{w}\bd_{X}=\bar{w}\circ\bd_{X}$.

Since $w[\Gamma']\subseteq F_{L_{h}M}$, we have that $(w(\varepsilon'(\psi)),w(\delta'(\psi)))\in\Omega^{L_{h}M}(F_{L_{h}M})$ for all $\psi\in\Gamma$ and $(\varepsilon',\delta')\in\tau'$. Since $\langle L_{h}M,F_{L_{h}M}\rangle$ is a reduced matrix, we have that
\[\begin{array}{rcl}
w(\varepsilon'(\psi))&=&w(\delta'(\psi))\\
\bar{w}(\varepsilon'(\psi))&=&\bar{w}(\delta'(\psi))\\
\bar{\rho_{M}}\circ\bar{w}(\varepsilon'(\psi))&=&\bar{\rho_{M}}\circ\bar{w}(\delta'(\psi))\\
\bar{v}\circ\bar{\rho_{X}}(\varepsilon'(\psi))&=&\bar{v}\circ\bar{\rho_{X}}(\delta'(\psi))
\end{array}\]

From $\bar{\rho}_{X}$ there is $\rho_{X}$ such that the square in the bellow diagram commutes, and then the diagram commutes:

\[\xymatrix{
F'(X)\ar[r]^{\rho_{X}}\ar[d]&F(X)\ar[r]^{v}\ar[d]&M\\
F'(X)/\Delta'\ar[r]_{\bar{\rho}_{X}}&F(X)/\Delta\ar[ur]_{\bar{v}}
}\]

With that we have $v\circ\rho_{X}(\varepsilon(\psi))=v\circ\rho_{X}(\delta(\psi))$, so $(v\circ\rho_{X}(\varepsilon(\psi)),v\circ\rho_{X}(\delta(\psi)))\in\Omega^{M}(F_{M})$ for all $\psi\in\Gamma'$. By algebraizability we have $v\circ\rho_{X}(\psi)\in F_{M}$ for all $\psi\in F_{M}$. By assumption $v\circ\rho_{X}(\varphi)\in F_{M}$. Thus $(v\circ\rho_{X}(\varepsilon(\varphi)),v\circ\rho_{X}(\delta(\varphi)))\in\Omega^{M}(F_{M})$. Since $\langle M,F_{M}\rangle$ is a reduced matrix, we have that
\[\begin{array}{rcl}
v\circ\rho_{X}(\varepsilon'(\varphi))&=&v\circ\rho_{X}(\delta'(\varphi))\\
\bar{v}\circ\bar{\rho}_{X}(\varepsilon'(\varphi)/\Delta')&=&\bar{v}\circ\bar{\rho}_{X}(\delta'(\varphi)/\Delta')\\
\bar{\rho}_{M}\circ\bar{w}(\varepsilon'(\varphi)/\Delta')&=&\bar{\rho}_{M}\circ\bar{w}(\delta'(\varphi)/\Delta')\\
\bd_{M}\circ\bar{\rho}_{M}\circ\bar{w}(\varepsilon'(\varphi)/\Delta')&=&\bd_{M}\circ\bar{\rho}_{M}\circ\bar{w}(\delta'(\varphi)/\Delta')\\
\bar{w}(\varepsilon'(\varphi)/\Delta')&=&\bar{w}(\delta'(\varphi)/\Delta')\\
w(\varepsilon'(\varphi))&=&w(\delta'(\varphi))
\end{array}\]

With that we have $(w(\varepsilon'(\varphi)),w(\delta'(\varphi)))\in\Omega^{L_{h}M}(F_{L_{h}M})$. Therefore $w(\varphi)\in F_{L_{h}M}$.

\vtres

{\bf (II)} One can use similar argument to prove the second part, i.e., given $a_1\in\sig(a)$, $\langle M_1,F_{M_1}\rangle\in Mod_{a}(a_1)=Matr^{*}_{a_1}$ and $\langle\Gamma'/\Delta',\varphi'/\Delta'\rangle\in Sen_{a'}(a')$ then:

\[\langle L_{h}(e_{1}^{*}(M_1)),F_{L_{h}(e_{1}^{*}(M_1))}\rangle\models'\langle\Gamma'/\Delta',\varphi'/\Delta'\rangle\ \ i\!f\!\!f\ \ \langle M_1,F_{M_1}\rangle\models_1\langle\check{e_{1}}\rho_{X}[\Gamma']/\Delta_{1},\check{f}\rho_{X}(\varphi')/\Delta_1\rangle\]

\qed

As a consequence of this theorem we have the abstract Glivenko's theorem between algebraizable logics.

\begin{Cor}\label{AbstGliInsAL}
For each Glivenko's context $\mathbb{G}=(h:a\to a',\rho)$, is associated an abstract Glivenko's theorem between $a$ and $a'$ i.e;  given $\Gamma'\cup\{\varphi'\}\subseteq F'(X)$ then
\[\rho_{X}[\Gamma']\vdash\rho_{X}(\varphi')\ \Leftrightarrow\ \Gamma'\vdash'\varphi'\]
\end{Cor}

\Dem

We know that for any algebraizable logic $a$,  $\vdash_a = \vdash_{Matr^{*}_{a}}$.
%CITAR TRECHO DO CAPITULO 1
 For any reduced matrix in $a'$ $\langle M',F_{M'}\rangle$ we have that $M'\in QV(a')$ and then there is $M\in QV(a)$ such that $L_{h}M \cong M'$ (see Remark \ref{GliDense}.(c)) , moreover $\langle L_{h}M,F_{L_{h}M}\rangle \cong \langle M',F_{M'}\rangle$. With that it is enough to prove that

\[\rho_{X}[\Gamma']\vdash_{Matr^{*}_{a}}\rho_{X}(\varphi')\ \Leftrightarrow\ \Gamma'\vdash_{Matr^{*}_{a'}}\varphi'\]

And that is equivalent to prove that for any $\langle M,F_{M}\rangle\in Matr^{*}_{a}$,

\[\langle M,F_{M}\rangle\models\langle\rho_{X}[\Gamma'],\rho_{X}(\varphi')\rangle\ \ i\!f\!\!f\ \ \langle L_{h}M,F_{L_{h}M}\rangle\models'\langle\Gamma',\varphi'\rangle\]

Or even,

\[\langle M,F_{M}\rangle\models\langle\rho_{X}[\Gamma']/\Delta,\rho_{X}(\varphi')/\Delta\rangle\ \ i\!f\!\!f\ \ \langle L_{h}M,F_{L_{h}M}\rangle\models'\langle\Gamma'/\Delta',\varphi'/\Delta'\rangle\]

But this last one follows from the previous theorem.\qed
\vtres

Now we present that the abstract Glivenko's theorem restricts to the classical Glivenko's theorem.

\begin{Ex}\label{traditionalGliv}
Let $\Sigma=(\Sigma_{n})_{n\in\omega}$ such that $\Sigma_{0}=\emptyset,\ \Sigma_{1}=\{\neg\},\ \Sigma_{2}=\{\longrightarrow\}$ and $\Sigma_{n}=\emptyset$ for all $n>2$. Let the map $h:IPC\to CPC$ such that $IPC$ and $CPC$ both are defined with the signature $\Sigma$,  $h(\neg)=\neg$ and $h(\longrightarrow)=\longrightarrow$, i.e., $h$ is the inclusion map from the intuitionistic propositional logic to the classical propositional logic. $IPC$ and $CPC$ are (Lindenbaum) algebraizable logics and $h$ is a morphism in $\cA_f$. Notice that $h^{*}$ is the identity functor and its restriction $h^{*}\!\!\rest:Bool\hookrightarrow Heyt$ has a left adjoint given by $L_{h}:Heyt\to Bool$ such that for any $A\in Heyt$, $L_{h}(A)=A_{\neg\neg}$ where is the boolean algebra of regular element, i.e., $a\in A$ such that $\neg\neg a=a$. The unit of this adjunction is $\partial_{A}:A\to A_{\neg\neg}$ such that $\partial_{A}(a)=\neg\neg a$ for all $A\in\Sigma-Str$. It is easy to see that this map define a natural transformation, moreover it has a natural transformation such that is a section given by $\rho_{A}:A_{\neg\neg}\to A$ where $\rho_{A}(a)= \neg\neg a = a$. Then we have that $(h:IPC\to CPC, \rho)$ is a Glivenko's context.

We know that $\psi\dashv_{CPC}\vdash\neg\neg \psi$ and then we have that $\psi/\Delta=\neg\neg \psi/\Delta$ where $\Delta=\{x\to y,\ y\to x\}$. Using the abstract Glivenko's theorem we have that given $\Gamma\cup\{\varphi\}$ set of formulas, then to prove that $\neg\neg \Gamma\vdash_{IPC}\neg\neg \varphi\ \Leftrightarrow\ \Gamma\vdash_{CPC}\varphi$ is enough to prove that for all matrix $\langle M,F_{M}\rangle\in Matr^{*}_{ICP}$,

\[\langle M_{\neg\neg},F_{M_{\neg\neg}}\rangle\models_{CPC}\langle\Gamma/\Delta,\varphi/\Delta\rangle\ \ i\!f\!\!f\ \ \langle M,F_{M}\rangle\models_{IPC}\langle\neg\neg\Gamma/\Delta,\neg\neg\varphi/\Delta\rangle\]

That is exactly the same to prove that

\[\langle L_{h}M,F_{L_{h}M}\rangle\models_{CPC}\langle\Gamma/\Delta,\varphi/\Delta\rangle\ \ i\!f\!\!f\ \ \langle M,F_{M}\rangle\models_{CPC}\langle\rho_{X}[\Gamma]/\Delta,\rho_{X}(\varphi)/\Delta\rangle.\]

This last follows from the previous corollary.

\end{Ex}

\begin{Obs} We believe that the notion of abstract Glivenko's theorem provided here, partially generalizes the approach that has been developed in \cite{To} In that paper, the author consider abstract Glivenko's theorem in the algebraizable logic setting (and also in some variants) but just relating logics defined over the {\em same signature} by means of an essentially idempotent formula with a free variable.
\end{Obs}

\subsection{The abstract Glivenko's theorem in \emph{InsLAL}}

We also have that a Glivenko's context induces an abstract Glivenko's theorem for $InsLAL$ and we present now.

In this subsection  we consider fixed: $a,a'$ Lindenbaum algebraizable logics, $G =(h:a\to a',\rho)$  a Glivenko's context
and a choice of isomorphisms \[\varepsilon_a:Obj(\sig'_{a}) \to \bigcup_{a_1 \in Obj(\sig'_a)} Hom_{\sig_a} (a, a_1)\]
given by $a_1 \mapsto ([e_1] : a \ra a_1)$

\begin{Cor}\label{LALprep}

For each $s = ([\varphi_0], \cdots, [\varphi_{n-1}], [\psi])$, a non-empty finite sequence in $F(\Sigma)/\dashv \vdash$ (the free $QV(a)$-structure on the set $X$) and each $(\tau, \Delta)$, an algebraizable pair of $a$,  where $\tau=\{(\varepsilon^{j},\delta^{j});\ j=1,...,m\ for\ some\ m\in\omega\}$, let
$$q(s, (\Delta, \tau)) :=  (([\varepsilon(\varphi_0)], [\delta(\varphi_0)]), \cdots, ([\varepsilon(\varphi_{n-1})],
[\delta(\varphi_{n-1})]), ([\varepsilon(\psi)], [\delta(\psi)]))$$
where the {\em notation} $([\varepsilon(\theta)], [\delta(\theta)])$ abbreviates  the pair of finite sequence of equivalence class of formulas: $([\varepsilon^j(\theta),  [\delta^j(\theta)])_j \ with \ j =1, \cdots, m.$.
Note that, as $a$ is a {\em congruential algebraizable  logic}, then: \\
(*) \  If $[\theta] = [\theta']$ (i.e., $\theta \dashv \vdash \theta'$), then $\delta(\theta) \dashv \vdash \delta(\theta')$ and  $\varepsilon(\theta) \dashv \vdash \varepsilon(\theta')$.Thus we have an well defined mapping $\varphi/\Delta \overset{t}\mapsto (\varepsilon(\varphi)/\Delta, \delta(\varphi)/\Delta)$ and $q(s, (\Delta, \tau))$ is well-defined;\\
(**) \ conversely, as $\varphi\dashv\vdash\Delta(\epsilon(\varphi),\delta(\varphi))$, then we have and well defined map $(\varepsilon(\varphi)/\Delta, \delta(\varphi)/\Delta) \overset{r}\mapsto \varphi/\Delta$ and $ r \circ t = id$.

Recall that $q_s :=  \{ q(s, (\tau, \Delta)) : (\tau, \Delta)$ is an algebraizable pair of $a_1\}$ and   $Sen'_{a}(a) := \{ q_s : s$  is a non-empty finite sequence in $F(\Sigma)\}$. Note that, by the above remark,  the mapping $s \overset{t} \mapsto q_s$ determine a {\em bijection} between the set of non-empty finite sequences in $F(\Sigma)/\Delta$ and $Sen'a(a)$

Then, in particular $\check{h}/\dashv \vdash :   F(\Sigma)/\dashv\vdash \to  F(\Sigma')/\dashv'\vdash$ has a section $\bar{\rho}_X :  F(\Sigma')/\dashv'\vdash \to  F(\Sigma)/\dashv\vdash$.  Let $s' = ([\varphi'_0], \cdots, [\varphi'_{n-1}], [\psi'])$ be a non-empty  finite sequence in $F(\Sigma')/\dashv'\vdash$ and $((\varepsilon',\delta'), \Delta')$ be an algebraizable pair of $a'$. Then  $\rho_*s' := ([\rho_X(\varphi'_0)], \cdots, [\rho_X(\varphi'_{n-1})], [\rho_X(\psi')])$ is a non-empty finite sequence in $F(\Sigma)/\dashv\vdash$ and the mapping
$q'_{s'} \in Sen'_{a'}(a')\  \mapsto \ q_{\rho_*s'} \in Sen'_{a}(a)$ is a section of the map on non-empty finite sequences induced by $\check{h}/\dashv\vdash:  F(\Sigma)/\dashv\vdash \ra F(\Sigma')/\dashv'\vdash$.

\end{Cor}

Now, we start providing the following

\begin{Prop}\label{Lind-Gliv-prop} Let $L_h : QV(a) \to QV(a')$ be the left adjoint of $h^\star\!\!\rest : QV(a') \to QV(a)$ as defined in Chapter 2 (see Proposition \ref{adjunct-prop}), then for each $M \in QV(a)$ and $q' \in Sent_a'(a')$,  the following  compatibility relation holds:
$$M\models^{a}\br q'\ \Leftrightarrow\ L_{h}M\models^{a'}q'$$

\end{Prop}

\Dem

$``\Rightarrow"$
Let $q'\in Sen'(a')$. Suppose that given $M\in QV(a)$, $M\models^{a}\br q'$. Given $w:X\to L_{h}(M)$ ($|L_{h}M|=|h^{\star}L_{h}M|$, we can consider $w:X\to h^{\star}L_{h}(M)$) such that
\[\bar{w}[\vepsilon'(\vphi'_{i})]=\bar{w}[\delta'(\vphi'_{i})],\ i=0,...,n-1\]

Look to diagram below:

\[\xymatrix{
&X\ar[dl]\ar@/_3pc/[ldd]_{v}\ar[dr]\ar@/^3pc/[rdd]^{w}&\\
\frac{FX}{\Delta}\ar@<1ex>[rr]^{\bd_{X}}\ar[d]^{\bar{v}}&&(\frac{F'X}{\Delta'})^{h}\ar[d]_{\bar{w}}\ar@<1ex>[ll]^{\br_{X}}\\
M\ar@<1ex>[rr]^{\bd_{M}}&&h^{\star}\!\!\rest L_{h}M\ar@<1ex>[ll]^{\br_{M}}
}\]

Consider $\bar{v}=\br_{M}\bar{w}\bd_{X}$ (there is $v:X\to M$ such that to be corresponding with $\bar{v}$). Hence $\bar{v}\br_{X}=\br_{M}\bar{w}\bd_{X}\br_{X}=\br_{M}\bar{w}$

\[\begin{array}{rcl}
\br_{M}\bar{w}[\vepsilon'(\vphi'_{i})]&=&\br_{M}\bar{w}[\delta'(\vphi'_{i})]\  (i=0,...,n-1)\\
\bar{v}\br_{X}[\vepsilon'(\vphi'_{i})]&=&\bar{v}\br_{X}[\delta'(\vphi'_{i})]\  (i=0,...,n-1)\\
\bar{v}\br_{X}[\vepsilon'(\vphi')]&=&\bar{v}\br_{X}[\delta'(\vphi')]\ Hypo.\\
\bd_{M}\bar{v}\br_{X}[\vepsilon'(\vphi')]&=&\bd_{M}\bar{v}\br_{X}[\delta'(\vphi')]\\
\bd_{M}\br_{M}\bar{w}[\vepsilon'(\vphi')]&=&\bd_{M}\br_{M}\bar{w}[\delta'(\vphi')]\\
\bar{w}[\vepsilon'(\vphi')]&=&\bar{w}[\delta'(\vphi')]
\end{array}\]

$w$ was taken arbitrary, so $L_{h}M\models^{a'}q$

$``\Leftarrow"$
Suppose that $L_{h}M\models^{a'}q'$. Let $v:X\to M$ such that $\bar{v}\br_{X}[\vepsilon'(\vphi'_{i})]=\bar{v}\br_{X}[\delta'(\vphi'_{i})]\ \forall\ i=0,...,n-1$.

Consider $\bar{w}=h^{\star}\!\!\rest L_{h}(\bar{v})$(exist $w:X\to h^{\star}\!\!\rest L_{h}M$ extends to $\bar{w}$). So $\br_{M}\bar{w}=\bar{v}\br_{X}$ and $\bar{w}\bd_{X}=\bd_{M}\bar{v}$. Therefore

\[\begin{array}{rcl}
\br_{M}\bar{w}[\vepsilon'(\vphi'_{i})]&=&\br_{M}\bar{w}[\delta'(\vphi'_{i})]\ i=0,...,n-1\\
\bd_{M}\br_{M}\bar{w}[\vepsilon'(\vphi'_{i})]&=&\bd_{M}\br_{M}\bar{w}[\delta'(\vphi'_{i})]\ i=0,...,n-1\\
\bar{w}[\vepsilon'(\vphi'_{i})]&=&\bar{w}[\delta'(\vphi'_{i})]\ i=0,...,n-1\\
\bar{w}[\vepsilon'(\vphi')]&=&\bar{w}[\delta'(\vphi')]\ Hypo.\\
\end{array}
\]

Hence $\bar{v}\br_{X}[\vepsilon'(\vphi')]=\br_{M}\bar{w}[\vepsilon'(\vphi')]=\br_{M}\bar{w}[\delta'(\vphi')]=\bar{v}\br_{X}[\delta'(\vphi')]$.

Then $M\models^{a}\br_{X} q'$
\qed

%%%%%%%%%%%%%%%%%%%%%%%%%%%%%%%%

We also have that a Glivenko's context induces an abstract Glivenko's theorem for $InsLAL$ and we present now and  Proposition \ref{Lind-Gliv-prop} above is part of it.

%From Proposition \ref{Lind-Gliv-prop} and Corollary \ref{AbstGliInsLAL} we are able to show

\begin{Teo}\label{MorGlivContInsLAL}
Let $a, a'$ be {\em Lindenbaum} algebraizable logics, then each $\mathbb{G}=(h:a\to a',\rho)$ \textbf{Glivenko's context} induces a institutions morphism $InsLAL_{a} \ra InsLAL_{a'}$. More precisely, fixing a choice of isomorphisms $\varepsilon: Obj(\sig'_{a}) \to Mor(\sig'_a)$, $a_1 \mapsto \varepsilon(a_1) = [e_1] : a \overset\cong\to a_1$,  we define a institution morphism $M_{(G,\varepsilon_a)} : InsLAL_{a}\ra InsLAL_{a'}$\footnote{Such induced morphisms are ``isomorphic", for different choices of isomorphisms $\varepsilon^0, \varepsilon^1$.}
\end{Teo}

\Dem

By simplicity, we will write $(G,\varepsilon$) for $(G,\varepsilon_{a})$. We will define
\[M_{(G,\varepsilon)}=\langle \Phi'^{(G,\varepsilon)},\alpha'^{(G,\varepsilon)},\beta'^{(G,\varepsilon)}\rangle\] (this will depend only on the choice of isomorphisms in the {\em domain} institution $InsLAL_{a}$):

$\Phi'^{(G,\varepsilon)} : \sig'_a \to \sig'_{a'}$: it is defined in the same way as $\Phi^{(G,\varepsilon)} : \sig_a \to \sig_{a'}$was defined in \ref{InsAL}.

%$ \bullet \ \Phi^\varepsilon_{h,\rho} : \sig_a \to \sig_{a'}$

%For each $[f] : a \to a$ (it is a isomorphism) and each $c'_n \in Sigma'_n$, let $\phi_n := \rho_X(c'_n(\vec{x})) \in F(\Sigma)[n]$, it is  such that $\vdash' \check{h}(\phi_n) \Delta' c'_n(\vec{x})$ (remember that $h$ is $\Delta$-dense). Define the $\cS_f$-morphism $f_{h} : \Sigma' \to \Sigma'$ as $f_h(c'_n) := \check{h}(\check{f}(\phi_n))$, then
%if we show that $f_h$ preserves $\vdash'$, then as $a'$ is congruential logic\\
%(I) $$[{f}_h] \circ[{h}] = [h] \circ [f]$$
% and, as $h$ is $\Delta$-dense, $[{f_h}]$ is the unique morphism satisfying (I). Moreover, as $[f]$ is an $Q(\cA^c_f)$-isomorphism it will follows that $[f_h]$ is a $Q(\cL_f)$-isomorphism and, as $h, f$ preserves algebraizing pairs, then (I) again ensures that $f_h$ is a $Q(\cA^c_f)$-isomorphism.

%Now suppose $\Gamma' \vdash' \psi'$, we must show that
%$\check`{f}_h[\Gamma'] \vdash' \check`{f}_h(\psi')$. As $\check{h}(rho_X(\phi) \dashv'\vdash \phi$, \\
%$\Gamma' \vdash' \psi' \ iff \ \check{h}[\rho_X[\Gamma']] \vdash' \check{h}(\rho_X(\psi'))$\\
%iff $\rho_{X}\partial_{X}[\rho_X[\Gamma']]\vdash\rho_{X}\partial_{X}(\rho_X(\psi'))$ \ (by Corollary \ref{AbstGliInsLAL}) \\
%iff $\check{f}\rho_{X}\partial_{X}[\rho_X[\Gamma']]\vdash\check{f} \rho_{X}\partial_{X}(\rho_X(\psi'))$ \ (because $[f]$ is an isomorphism) \\
%iff $\check{f}[\rho_X[\Gamma']]\vdash\check{f} (\rho_X(\psi'))$ \ (because $\check{h}[\rho_X (\theta) \dashv' \vdash \theta$ ) \\

\vtres

Now the definition of $\alpha'^{(G,\varepsilon)}$.

Firstly for $a$ we have $\alpha'^{(G,\varepsilon)}(a):Sen_{a'}\circ\Phi'^{(G,\varepsilon)}(a)=Sen_{a'}(a')\to Sen_{a}(a)$ is the mapping  $q'_{s'} \in Sen'_{a'}(a')\  \mapsto \ q_{\rho_*s'} \in Sen'_{a}(a)$, as defined in \ref{LALprep}.

%Consider $q'=(([\varepsilon'(\varphi'_{0})],[\delta'(\varphi'_{0})]),...,([\varepsilon'(\varphi'_{n-1})],[\delta'(\varphi'_{n-1})]),([\varepsilon(\varphi)],[\delta(\varphi')]))$. Due to \ref{GliDense} and since $a$ and $a'$ are algebraizable logics, there is $(\Delta,(\varepsilon,\delta))$ algebraizable pair in $a$ and there is $\{\varphi_{0},...,\varphi_{n-1},\varphi\}\subseteq FX$ such that $\check{h}(\varepsilon(\varphi_{i}))\dashv'\vdash\varepsilon'(\varphi'_{i})\ \forall\ i=0,...,n-1$ and $\check{h}(\varepsilon(\varphi))\dashv'\vdash\varepsilon'(\varphi)$. The same conditions for $\delta$ and $\delta'$. Therefore
%\[q'=(([\check{h}(\varepsilon(\varphi_0))],[\check{h}(\delta(\varphi_0))]),...,([\check{h}(\varepsilon(\varphi_{n-1}))],[\check{h}(\delta(\varp%hi_{n-1}))]),([\check{h}(\varepsilon(\varphi)),\check{h}(\delta(\varphi))]))\]
%Moreover
%\[q'=((\bar{\partial}_{X}[\varepsilon(\varphi_0)],\bar{\partial}_{X}[\delta(\varphi_0)]),...,(\bd_{X}[\vepsilon(\vphi_{n-1})],\bd_{X}[\delta(\vphi_{n-1})]),(\bd_{X}[\vepsilon(\vphi)],\bd_{X}[\delta(\vphi)]))\]

%Defines
%\[\alpha_h(a)(q')=\br_{X}q'=((\br_{X}[\varepsilon'(\varphi'_{0})],\br_{X}[\delta'(\varphi'_{0})]),...,(\br_{X}[\varepsilon'(\varphi'_{n-1})],\br_{X}[\delta'(\varphi'_{n-1})]),(\br_{X}[\varepsilon(\varphi)],\br_{X}[\delta(\varphi')]))\]

For an arbitrary $a_{1}\in\sig'_{a}$ we define $\alpha'^{(G,\varepsilon)}(a_{1}):Sen_{a'}(a')\to Sen_{a}(a_1)$ by for $q'\in Sen_{a}(a')$,
\[\alpha'^{(G,\varepsilon)}(a')(q')= \br_X^{a_1}(q')\]

such that for each component of $\br_{X}^{a_1}q'$ is $([\check{e_1}\rho_{X}(\varepsilon'(\varphi_{k}))],[\check{e_1}\rho_{X}(\delta'(\varphi_{k}))])$ for $k=1,...,n-1$ and the last component is $([\check{e_1}\rho_{X}(\varepsilon'(\varphi))],[\check{e_1}\rho_{X}(\delta'(\varphi))])$.
This defines a natural transformation. Indeed, first observe that the diagram below commutes:

\[\xymatrix{
F(\Sigma')/\Delta'\ar[r]^{\br_X}\ar[d]_{[\check{g}^{g}_{\varepsilon}]}&F(\Sigma)/\Delta\ar[r]^{[\check{e}_{1}]}\ar[d]^{[\check{g}_{\varepsilon}]}&F(\Sigma_1)/\Delta_1\ar[d]^{[\check{g}]}\\
F(\Sigma')/\Delta'\ar[r]_{\br_X}&F(\Sigma)/\Delta\ar[r]_{[\check{e}_2]}&F(\Sigma_2)/\Delta_2
}\]

then we have the following diagram commuting:

\[\xymatrix{
Sen_{a'}(\Phi'(a_1))\ar[d]_{[\check{g}_{\varepsilon}^{h}]}\ar[rr]^{\alpha'^{(G,\varepsilon)}(a_1)}&&Sen_{a}(a_1)\ar[d]^{[\check{g}]}\\
Sen_{a'}(\Phi'(a_2))\ar[rr]_{\alpha'^{(G,\varepsilon)}(a_2)}&&Sen_{a}(a_2)
}\]

\vtres

Let now to define $\beta'^{(G,\varepsilon)}$. For $a$ we define  $\beta'^{(G,\varepsilon)}:Mod'_{a}\Rightarrow Mod'_{a'}\circ(\Phi'^{(G,\varepsilon)})^{op}$ is define as:

$\beta'^{(G,\varepsilon)}(a)=L_{h}:QV(a)=Mod'_{a}(a)\to Mod'_{a'}(\Phi'^{(G,\varepsilon)}(a))=QV(a')$

The corresponding definition works for an arbitrary $a_1\in\sig_{a}$ because since $a$ and $a_1$ are $Q^{c}_{f}$-isomorphic, we have by \ref{QV-iso-prop} that $QV(a)$ and $QV(a_1)$ are isomorphic. I.e., $\beta'^{(G,\varepsilon)}(a_1)=L_{h_1}:QV(a_1)=Mod'_{a}(a_1)\to Mod'_{a'}(\Phi'^{(G,\varepsilon)}(a_1))=QV(a')$, where $(a \overset{[h]}\ra a')\ =\ (a \overset{[e_1]}\ra a_1 \overset{[h_1]}\ra a')$.  This defines a natural transformation. Indeed, notice that the following diagram commutes:

\[\xymatrix{
QV(a_1)&QV(a')\ar[l]_{h_{1}^{\star}\!\!\rest}\\
QV(a_2)\ar[u]^{g^{\star}\!\!\rest}&QV(a')\ar[l]^{h_{2}^{\star}\!\!\rest}\ar[u]_{(g_{\varepsilon}^{h})^{\star}\!\!\rest}
}\]

And then we have the following diagram commuting:

\[
\xymatrix{
QV(a_1)\ar[r]^{\beta_{a_1}}&QV(a')\\
QV(a_2)\ar[r]_{\beta_{a_2}}\ar[u]^{g^{\star}\!\!\rest}&QV(a')\ar[u]_{(g^{h}_{\varepsilon})^{\star}\!\!\rest}
}\]

\vtres

On the compatibility condition. First for the logic $a$ we must guarantee that $M\models^{a}\br q'\ \Leftrightarrow\ L_{h}M\models^{a'}q'$: this is the content of Proposition \ref{Lind-Gliv-prop}.

For an arbitrary logic $a_1\in\sig'_a$ we must to prove that for any $M_1\in QV(a_1)$ and $q_{s'}\in Sen(a')$:

\[\beta'^{(G,\varepsilon)}(a_1)(M_1)\models^{a'}q_{s'}\ \Leftrightarrow\ M_1\models^{a_1}\alpha'(G,\varepsilon)(a_1)(q_{s'})\]
in other notation
\[L_{h_1}(M_1)\models^{a'}q_{s'}\ \Leftrightarrow\ M_1\models^{a_1}\br^{a_1}_{X}(q_{s'}).\]

In fact, since $[\Phi(e_{1})]$ is an isomorphism, we have that $\Phi(e_{1})^{\star}\!\!\rest$ is an isomorphism. Therefore:

\[\begin{array}{rcl}
L_{h_{1}}(M_1)\models^{a'}q_{s'}&\Leftrightarrow&\Phi(e_{1})^{\star}L_{h_1}(M_1)\models q_{(\Phi(e_{1})^{+})^{-1}(s_1)}\\
&\Leftrightarrow&L_{h}(e_{1}^{\star}(M_1))\models q_{(\Phi(e_{1})^{+})^{-1}(s_1)}\\
&\Leftrightarrow&e_{1}^{\star}(M_1)\models\alpha'^{(G,\varepsilon)}(a)(q_{(\Phi(e_{1})^{+})^{-1}(s_1)})\\
&\Leftrightarrow&M_{1}\models^{a_1} e_{1}^{+}\alpha'^{(G,\varepsilon)}(a)(q_{(\Phi(e_{1})^{+})^{-1}(s_1)})\\
&\Leftrightarrow&M_{1}\models^{a_1}\alpha_{a_1}(q_{s'}).
\end{array}\]

\qed

\begin{Cor}\label{AbstGliInsLAL}
For each Glivenko's context $\mathbb{G}=(h:a\to a',\bar{\rho})$, is associated an abstract Glivenko's theorem between $a$ and $a'$ i.e;  given $\Gamma'\cup\{\varphi'\}\subseteq F'(X)$ then
\[\rho_{X}[\Gamma']\vdash\rho_{X}(\varphi')\ \Leftrightarrow\ \Gamma'\vdash'\varphi'\]
\end{Cor}

\Dem

Firstly, remark that it is enough consider $\Gamma$ finite. Because $a$ and $a'$ are algebraizable logics, and  $h$ preserves algebraizing pairs, it is enough to show that
\[
\begin{array}{rcl}
\{\vepsilon(\rho_{X}(\psi'))\approx\delta(\rho_{X}(\psi')),\psi'\in\Gamma'\}&\models_{QV(a)}&\vepsilon(\rho_{X}(\vphi'))\approx\delta(\rho_{X}(\vphi'))\\
&\Updownarrow&\\
\{\vepsilon'(\psi')\approx\ \delta'(\psi'),\psi'\in\Gamma'\}&\models_{QV(a')}&\vepsilon'(\vphi')\approx \delta'(\vphi')\\
%&\Updownarrow &\\
%\{\check{h}(\vepsilon(\psi))\approx\check{h}(\delta(\psi)),\psi\in\Gamma\}&\models_{QV(a')}&\check{h}(\vepsilon(\vphi))\approx\check{h}(\delta(\vphi))
\end{array}\]

%Due to $a$ and $a'$ are Lindenbaum algebraizable,

%Thus, is enough to show that

%\[\begin{array}{rcl}
%\{\br_{X}\bd_{X}[\vepsilon(\psi)]\approx\br_{X}\bd_{X}[\delta(\psi)],\psi\in\Gamma\}&\models_{QV(a)}&\br_{X}\bd_{X}[\vepsilon(\vphi)]\approx\br_{X}\bd_{X}[\delta(\vphi)]\\
%&\Updownarrow&\\
%\{[\check{h}(\vepsilon(\psi))]\approx[\check{h}(\delta(\psi))],\psi\in\Gamma\}&\models_{QV(a')}&[\check{h}(\vepsilon(\vphi))]\approx[\check{h}(%\delta(\vphi))]
%\end{array}\]

Consider $\Gamma=\{\psi_{0},...,\psi_{n-1}\}$,
$s' =  (\psi'_0/\Delta', ...,\psi'_{n-1}/\Delta', \varphi'/\Delta')$. Then:\\
(i) $q'= q'_{s'}$ is determined by any of its elements \[(([\vepsilon'(\psi'_{0})],[\delta'(\psi'_{0})]),...,([\vepsilon'(\psi'_{n-1})],[\delta'(\psi'_{n-1})])([\vepsilon'(\vphi')],[\delta'(\vphi)]));\]

(ii) $\alpha(a)(q') = q_{\rho_*s'}$ is determined by any of its elements \[(([\vepsilon(\rho_X(\psi'_{0}))],[\delta(\rho_X(\psi'_{0}))]),...,([\vepsilon(\rho_X(\psi'_{n-1}))],[\delta(\rho_X(\psi'_{n-1}))])([\vepsilon(\rho_X(\vphi'))],[\delta(\rho_X(\vphi))]))\]

Thus we have to show:

 \[(\forall M \in QV(a), M\models^{a}\br_{X}q') \ \Leftrightarrow\ (\forall M' \in QV(a') M'\models^{a'}q')\]

By Remark \ref{GliDense}.(c) given $M'\in QV(a')$ there is $M\in QV(a)$ such that
$L_{h}(M)\cong M'$.

With this, it is enough to show that for every $M\in QV(a)$, \[M\models^{a}\br_{X}q'\Leftrightarrow L_{h}M\models^{a'}q'\]

And this last equivalence is established the Proposition \ref{Lind-Gliv-prop} above.

\qed

\begin{Obs}
Since the CPC and IPC are Lindenbaum algebraizable logic, one can see that the example \ref{traditionalGliv} follows a consequence of the abstract Glivenko's theorem fo $InsLAL$ as well as the abstract Glivenko's theorem for $InsAL$.
\end{Obs}
%ENTAO, EU RESOLVI ESCREVER O COROLARIO DESSA FORMA PQ A FORMA ANTERIOR TINHA UM PROBLEMA. O $\rho_{X}$ NAO PODE SER APLICADO A FORMULAS DA LOGICA $a$ E SIM DA LOGICA $a'$. MAS ACHO QUE CONTINUA FAZENDO SENTIDO SE PENSARMOS NAS LOGICAS INTUICIONISTA E CLASSICA. LA, PELO Q ESTOU PENSANDO, $\partial_{X}$  EXATAMENTE A DUPLA NEGACAO E $\rho_{X}$ ACHO Q  A APLICACAO IDENTIDADE. NAO SEI... ACHO QUE  ISSO. CONVERSAMOS NA SEGUNDA.

%%%%%%%%%%%%%%%%%%%%%%%%%%%%%%%%%%%%%%%%%

\begin{Obs} A simple analysis of the derivations of ``logical" forms of Glivenko's Theorem (Corolaries \ref{AbstGliInsAL} and \ref{AbstGliInsLAL}) from the corresponding ``instituitional" form of Glivenko's Theorem  (Theorems \ref{MorGlivContInsAL} and \ref{MorGlivContInsLAL}), i.e. the existence of certain (induced) morphisms of institutions make clear that the latter form is stronger than the former one. We can interpret this as another evidence\footnote{Beside the nice approach of  the identity problem for (algebraizable) propositional logics: ``a logic is an institution, thus manifested through many signatures".} of the (virtually unexplored)  relevance of institution theory in propositional logic.
\end{Obs}

\section{Category of algebraizable logics with Glivenko's morphisms}

In this section we present that the definition of Glivenko's context given in \ref{GlivenkoContext} offer more information about the relationship of logics, it give us a category of algebraizable logics such that the morphisms are Glivenko's contexts, i.e., the objects are the same of in $\cA_f$ and given $a$ and $a'$ algebraizable logics, a Glivenko's morphism is a Glivenko's context $(h:a\to a',\rho)$. Denote by $G\cA_f$ this category.

\begin{Teo}
$G\cA_f$ is a category
\end{Teo}

\Dem

In this category the composition is the usual, i.e., given $G=(h:a\to a',\rho)$ and $G'=(h':a'\to a'',\rho')$, we have that $G'\circ G=(h'\circ h:a\to a'',\rho'\bullet\rho)$ where $(\rho'\bullet\rho)_{M}=\rho_{M}\circ\rho'_{L_{h}M}$ (this is natural in $M \in QV(a)$). In order to prove that the composition is well defined, we must to prove that $\rho'\circ\rho$ a section for the unit of the adjunction $L_{h'\circ h}\dashv (h'\circ h)^{*}$. The composition of adjunctions is a adjunction and $\partial'\circ\partial$ is its the unit. Remember that $L_{h'\circ h}=L_{h'}\circ L_{h}$ (an strict equality, with the choice of adjoints given in Chapter 2, as quotients) and $(h'\circ h)^{*}=h^{*}\circ h'^{*}$. Then we have that $M\overset{(\partial'\circ\partial)_{M}}\longrightarrow h^{*}h'^{*}L_{h'}L_{h}M =M\overset{h^{*}(\partial'_{L_{h}M})\circ\partial_{M}}\longrightarrow h^{*}h'^{*}L_{h'}L_{h}M=M\overset{\partial'_{L_{h}M}\circ\partial_{M}}\longrightarrow h^{*}h'^{*}L_{h'}L_{h}M$. Then we have
\[\begin{array}{rcl}
(\partial'\circ\partial)_{M}\circ(\rho'\bullet\rho)_{M}&=&(\partial'_{L_{h}M}\circ\partial_{M})\circ(\rho_{M}\circ\rho'_{L_{h}M})\\
&=&\partial'_{L_{h}M}\circ(\partial_{M}\circ\rho_{M})\circ\rho'_{L_{h}M}\\
&=&\partial'_{L_{h}M}\circ\rho'_{L_{h}M}\\
&=&Id_{L_{h}M}.
\end{array}\]

Thus $(\rho'\bullet\rho)_{M}$ is a section for $(\partial'\circ\partial)_{M}$ for all $M\in\Sigma-Str$. Clearly there is the identity Glivenko's context for an algebraizable logic $a$ given by $(Id_{a}:a\to a,\rho=(Id_{M})_{M\in\Sigma-Str})$. To prove the associativity let $G=(h:a\to a',\rho),\ G'=(h':a'\to a'',\rho')$ and $G''=(h'':a''\to a''',\rho'')$ be Glivenko's morphisms (Glivenko's context). Since $\cA_f$ is a category we have that $h''\circ(h'\circ h)=(h''\circ h')\circ h$. Remains to prove that $\rho''\bullet(\rho'\bullet\rho)=(\rho''\bullet\rho')\bullet\rho$. Let $M\in\Sigma-Str$, then

\[\begin{array}{rcl}
(\rho''\bullet(\rho'\bullet\rho))_{M}&=&(\rho'\bullet\rho)_{M}\circ\rho''_{L_{h'\circ h}M}\\
&=&(\rho_{M}\circ\rho'_{L_{h}M})\circ\rho''_{L_{h'\circ h}M}\\
&=&(\rho_{M}\circ\rho'_{L_{h}M})\circ\rho''_{L_{h'}\circ L_{h}M}\\
&=&\rho_{M}\circ(\rho'_{L_{h}M}\circ\rho''_{L_{h'}L_{h}}M)\\
&=&\rho_{M}\circ(\rho''\bullet\rho')_{L_{h}M}\\
&=&((\rho''\bullet\rho')\bullet\rho)_{M}
\end{array}\]

Therefore $G\cA_f$ is a category\qed
\vtres

The theorems  \ref{MorGlivContInsAL} and \ref{MorGlivContInsLAL}  say that for any Glivenko's context there is a institution morphism associated, more precisely, given a Glivenko's context $(h:a\to a',\rho)$ and a choice of isomorphisms $\varepsilon_a:Obj(\sig_{a}) \to \bigcup_{a_1 \in Obj(\sig_a)} Hom_{\sig_a} (a, a_1)$, we have a institution morphism $\langle \Phi_{G, \varepsilon},\alpha_{G,\varepsilon},\beta_{G,\varepsilon}\rangle$. Notice that there are more than one possible choice for the family $(\varepsilon_a)_{a \in |\cA_f|}$, but the application below still define a functor.

%%%%%%%COMPLETAR COM A PROVA!!!%%%%%%%

 %$\bullet$ identityIf $G = (id_a, id)$ and $\varepsilon:Obj(\sig_{a}) \to \bigcup_{a_1 \in Obj(\sig_a)} Hom_{\sig_a} (a, a_1)$ is a choice of isos:\\
%- $\Phi_{G, \varepsilon} = Id_{\sig_a}$;\\
%- $\alpha_{G,\varepsilon} =

%compsotion

\[\begin{array}{rrcl}
\mathcal{G}_{\varepsilon}:&G\cA_f&\to&{\bf Inst}\\
&a&&InsAL_{a}\\
&(h,\rho)\downarrow&\mapsto&\downarrow\langle\Phi_{(G,\varepsilon)},\alpha_{(G,\varepsilon)},\beta_{(G,\varepsilon)}\rangle\\
&a'&&InsAL_{a'}
\end{array}\]

\vtres

Another natural functor that arise is $\mathcal{U}:G\cA_f\to \cA_f$ such that $\mathcal{U}((h:a\to a',\rho))= (h:a\to a')$ for any Glivenko's context $(h:a\to a',\rho)$.

Naturally, we can defined in analogous way a (full) subcategory $G\cA^c_f \sub G\cA_f$, with objects being the Lindenbaum algebraizable logics and, for each choice of isomorphisms $(\varepsilon_a)_{a \in |\sig'_{a}|}$, we get a functor:

\[\begin{array}{rrcl}
\mathcal{G}^c_{\varepsilon}:&G\cA^c_f&\to&{\bf Inst}\\
&a&&InsLAL_{a}\\
&(h,\rho)\downarrow&\mapsto&\downarrow\langle\Phi'_{(G,\varepsilon)},\alpha'_{(G,\varepsilon)},\beta'_{(G,\varepsilon)}\rangle\\
&a'&&InsLAL_{a'}
\end{array}\]

%CUIDADO DINATURAL? COMORFISMO E MORFISMO!
%Moreover, the family of comorphism of institutions described in Remark \ref{compar- inst}, defines a natural transformation between the parallel functors $\eta_{\varepsilon}: \mathcal{G}^c_{\varepsilon\rest}  \Rightarrow \mathcal{G}_{\varepsilon} \circ incl$, where $\varepsilon\rest$ is the restriction to the Lindenbaum algebraizable logics, of (family of) choice of isomorphisms.

%%%%%%%COMPLETAR ISTO TAMBEM COM A PROVA!!!%%%%%%%

Once established those relations we have the following diagram that represents the relation among the categories studied in this thesis.

%\[
%\xymatrix{
%\pi-Inst\ar@/_1pc/[rr]&&Inst\ar@/_1pc/[ll]&&\mathcal{F}i\ar@/^1pc/[dd]\\
%&&G\cA_f\ar[u]\ar[d]&&\\
%Lind(\cA_f)\ar@{>->}[rr]\ar[uu]&&\cA_f\ar@{>->}[rr]\ar[uull]&&\Lf\ar@/^1pc/[uu]\ar[uullll]
%}\]

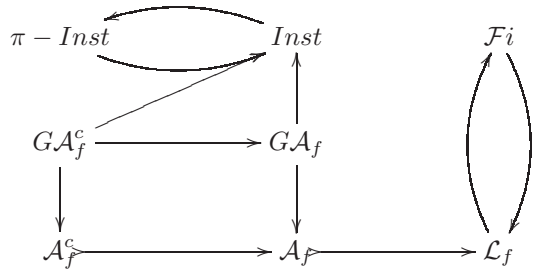
\[
\xymatrix{
\pi-Inst\ar@/_1pc/[rr]&&Inst\ar@/_1pc/[ll]&&\mathcal{F}i\ar@/^1pc/[dd]\\
G\cA^{c}_{f}\ar[urr]\ar[d]\ar[rr]&&G\cA_f\ar[u]\ar[d]&&\\
\cA^{c}_{f}\ar@{>->}[rr]&&\cA_f\ar@{>->}[rr]&&\Lf\ar@/^1pc/[uu]
}\]

On the other hand,we saw that the categories $\cL_s$ and $\cL_f$ determines institutions and $\pi$-institutions. Having in mind the adjunctions $\cL_s \rightleftarrows \cL_f  \rightleftarrows \cF_i$, we believe that is possible establish a (extended) direct relation from  $\mathcal{F}i$ to ${\bf Inst}$ and $\Pi-{\bf Inst}$. This is part of the future works on this paper.

%\end{thebibliography}

\end{document}